\documentclass[a4paper,11pt]{article}
\usepackage[plainpages=false]{hyperref}
\usepackage{amsfonts,latexsym,rawfonts,amsmath,amssymb,amsthm,mathrsfs}
\usepackage{amsmath,amssymb,amsfonts,latexsym,lscape,rawfonts}

\usepackage[all]{xy}
\usepackage{eufrak}
\usepackage{makeidx}         
\usepackage{graphicx,psfrag}

\usepackage{array,tabularx}

\usepackage{setspace}

\newtheorem{thm}{Theorem}[section]
\newtheorem{cor}[thm]{Corollary}
\newtheorem{lem}[thm]{Lemma}

\newtheorem{prop}[thm]{Proposition}
\newtheorem{defi}[thm]{Defintion}
\newtheorem{prob}[thm]{Problem}

\theoremstyle{remark}
\newtheorem{rmk}[thm]{Remark}

\theoremstyle{definition}

\def \C {\mathbb C}

\def \R {\mathbb R}
\def \D {\mathbb D}

\def \P {\mathbb P}
\def \A {\mathcal A}
\def \U {\mathcal U}
\def \V {\mathcal V}

\def \M {\mathcal M}
\def \W {\mathcal W}

\def \NB {\mathcal N}

\def \D {\mathcal D}
\def \J {\mathcal J}
\def \H {\mathcal H}
\def \G {\mathcal G}
\def \L {\mathcal L}
\def \g {\mathfrak g}

\def \m {\mathfrak m}

\def \p {\partial}
\def \pb {\bar{\partial}}
\def  \OT   {\Omega^{0,1}(T^{1,0})}
\def  \OTS   {\Omega_S^{0,1}(T^{1,0})}

\def  \OOTS  {\Omega_S^{0,2}(T^{1,0})}
\def \tgamma {\widetilde{\gamma}}

\begin{document}
\title{Calabi flow, Geodesic rays, and uniqueness of constant scalar curvature K\"ahler metrics}
\author{Xiuxiong Chen and Song Sun}
\date{}

\maketitle
\begin{abstract} We prove that constant scalar curvature K\"ahler metric
``adjacent" to a fixed K\"ahler class is unique up to isomorphism.
This extends the uniqueness theorem of  Donaldson and Chen-Tian,
and formally fits into the infinite dimensional G.I.T picture
described by Donaldson. We prove that the Calabi flow near a cscK
metric exists globally and converges uniformly to a cscK metric in
a polynomial rate. Viewed in fixed a K\"ahler class, the Calabi
flow is also shown to be asymptotic to a smooth geodesic ray at
infinity. This latter fact is also interesting in the finite
dimensional case, where we show that the downward gradient flow of
the Kempf-Ness function in a semi-stable orbit is asymptotic to
the direction of optimal degeneration.
\end{abstract}

\tableofcontents
\section{Introduction}
The Kempf-Ness theorem relates complex quotient to symplectic
reduction. Suppose a compact connected group $G$ acts on a compact
K\"ahler manifold $X$. We assume the action preserves the K\"ahler
structure, with a moment map $\mu: X\rightarrow \g^*$. Then the
action extends to a holomorphic action of the complexified group
$G^{\C}$. Under proper hypothesis the notion of stability could be
defined. Then the Kempf-Ness theorem says that as sets:
$$X^{ss}/G^{\C}\simeq\mu^{-1}(0)//G.$$ To be more precise,\\

(1). A $G^{\C}$-orbit is poly-stable if and only if it contains a
zero of the moment map. The zeroes within it form a unique $G$ orbit. \\

(2). A $G^{\C}$-orbit is semi-stable if and only if its closure
contains a zero of the moment map. We call such a zero  a
\emph{de-stabilizer} of the original $G^{\C}$ orbit. The
de-stabilizers all lie in the unique poly-stable orbit in the
closure of the
original orbit.\\

In K\"ahler geometry according to S. Donaldson(\cite{D1})(see also
\cite{Fu}) the problem of finding cscK (constant scalar curvature
K\"ahler) metrics formally fits  into a
 similar picture. However the spaces involved are infinite dimensional.
Given a compact K\"ahler manifold
 $(M, \omega, J)$, denote by
  $\G$ the group of Hamiltonian diffeomorphisms of $(M,
 \omega)$ and by $\J$ the space consists of almost complex structures on $M$
which are compatible with $\omega$.
 $\J$ admits a natural K\"ahler structure which is invariant under the
action of $\G$. The moment map is given by the Hermitian scalar
  curvature. The complexification of $\G$ may not exist, since $\G$ is
  infinite dimensional. Nevertheless, it still makes sense talking
  about the orbits of $\G^{\C}$--it is simply the leaf of the
  foliation obtained by complexifying the infinitesimal actions of $\G$.
  Then the $\G^{\C}$ leaf of an integrable complex structure can
  be viewed as a principal $\G$-bundle over the K\"ahler class
  $[\omega]$.
  Thus an analogue to the Kempf-Ness theorem should relate the stability of
the leaves to the existence of cscK
  metrics in the corresponding K\"ahler class. This was made more precise as
the \emph{Yau-Tian-Donaldson conjecture}(see
   \cite{Th}).
   The notion of ``stability" in this case is  the so-called ``K-stability", see \cite{Ti1}, \cite{D5}.
   There are also other related notion of stability, see for example \cite{RT}, \cite{Pa}, etc. \\

  Note that the Kempf-Ness theorem consists of both the existence and
  uniqueness part. It is known that the existence of cscK metrics implies
various kinds of stability, however the
  converse is fairly difficult, due to
  the appearance of fourth order non-linear P.D.E's. Recently
  Donaldson(\cite{D6}) proved a  general result that the conjecture is true
  for toric surfaces. The uniqueness part corresponding to the
  poly-stable case is known by\\

\begin{thm}(Donaldson\cite{D3}, Chen-Tian\cite{CT})\label{CT} Constant
scalar K\"ahler metric in a fixed
K\"ahler class, if exists, is unique up to holomorphic isometry.\\
\end{thm}
\begin{rmk} When the manifold is Fano, the uniqueness of
K\"ahler-Einstein metrics was previously proved by
Bando-Mabuchi(\cite{BM}), and  it was later generalized to the
case of K\"ahler-Ricci solitons by Tian-Zhu(\cite{TZ1}). The
uniqueness of cscK metrics was first proved by the first author in
the case when $c_1(X)\leq0$(\cite{Ch1}).
\end{rmk}

The purpose of this paper is to prove the uniqueness in the
semi-stable case.

\begin{thm} \label{Theorem 1} If there are two cscK structures $J_1$ and
$J_2$ both lying in the ($C^{\infty}$) closure of the $\G^{\C}$
leaf of a complex structure $J\in \J^{int}$, then there is a
symplectic diffeomorphism $f$ such that $f^*J_1=J_2$.
\end{thm}

\begin{defi} Let $(M, \omega, J)$ be a K\"ahler manifold and $\H$ be the
space of K\"ahler metrics in the K\"ahler class of $\omega$. We
say another K\"ahler structure $(\omega', J')$ on $M$ is adjacent
to $\H$ if there is a sequence of K\"ahler metrics $\omega_i\in
\H$ and diffeomorphisms $f_i$ of $M$ such that
$$f_i^*\omega_i\rightarrow\omega', f_i^*J\rightarrow J'$$ in
$C^{\infty}$ sense. So in particular, the corresponding sequence
of Riemannian metrics $g_i$ converges to $g'$ in the
Cheeger-Gromov sense. Similarly, let $(M, J)$ be a Fano manifold.
We say another complex structure $J'$ on $M$ is adjacent to $J$ if
there is a sequence of diffeomorphisms $f_i$ such that
$$f_i^*J\rightarrow J'.$$
\end{defi}

\begin{rmk} The above definition is related to the ``jumping"
phenomenon of complex structures, i.e. the space of isomorphism
classes of complex structures on a fixed manifold is in general
not Hausdorff. As a simple example, we can consider the blown-up
of $\P^2$ at three points $p_1$, $p_2$, and $p_3$. The underlying
differential manifold is fixed, and a choice of the three points
defines a complex structure. A choice of three points in a general position
gives rise to the same complex structure, while a choice of three
points on a line provides an example of an adjacent complex structure.
\end{rmk}

It follows theorem \ref{Theorem 1} that

\begin{thm} \label{Theorem 2}Let $(M,\omega,J)$ be a K\"ahler manifold.
Assume $[\omega]$ is integral. Suppose there are two csc K\"ahler
structures $(\omega_1, J_1)$ and $(\omega_2, J_2)$ both adjacent
to the K\"ahler class of $(\omega, J)$, then they are isomorphic.
\end{thm}

\begin{cor} \label{Corollary 1}Let $(M, J)$ be a Fano manifold. Suppose
there are two complex structures $J_1$ and $J_2$ both adjacent to
$J$ and both admitting K\"ahler-Einstein metrics, then $(M, J_1)$
and $(M, J_2)$ are bi-holomorphic.
\end{cor}

\begin{rmk} After finishing this paper, we learned that our theorem
\ref{Theorem 2} and corollary \ref{Corollary 1} partially
confirmed a conjecture of G. Tian(\cite{Ti2}) in the case of
constant scalar curvature K\"ahler metric.
\end{rmk}

The main technical ingredient in the proof of the above theorems
is to obtain some $C^0$ bound. We shall study the asymptotic
behavior of the Calabi flow near a cscK metric. The global
existence and convergence are established by using the Lojasiewicz
inequality which controls the gradient of a real analytic function
near a critical point. Suppose now we have two cscK metrics
adjacent to a fixed K\"ahler class, then there are two Calabi
flows in the neighborhoods of the corresponding cscK metrics.
Since the Calabi flow decreases geodesic distance, we get a bound
on the two Calabi flows in terms of geodesic distance. It is not
known whether this bound implies $C^0$ bound automatically. Here
we get around this difficulty by showing that the previous Calabi
flow is \emph{asymptotic} to a smooth geodesic ray. This involves
a local study of the infinite dimensional Hamiltonian action of
$\G$, which is the main technical part of this paper. We shall
first look at the analogous finite dimensional problem. Finally we
are able to derive $C^0$ bound for the two
parallel geodesic rays.\\

 The organization of this paper
is as follows. In section 2, we review Donaldson's infinite
dimensional moment map picture in K\"ahler geometry, and recall
some known results for our later use. In section 3, we state the
Lojasiewicz inequality and ``Lojasiewicz arguments" for the
gradient flow of a real analytic function. In section 4, we prove
that in the finite dimensional case, the Kempf-Ness flow for a
semi-stable point is asymptotic to a rational geodesic ray. In
section 5, we study the stability of the Calabi flow near a cscK
metric when the complex structure is deformed. In section 6, we
generalize the arguments in section 4 to the infinite dimensional
setting by considering the ``reduced" Calabi flow. In section 7,
the relative $C^0$ bound for two smooth parallel geodesic rays
tamed by bounded geometry is derived. In section 8, we prove the
main theorems. In Section 9, we shall discuss some further
problems related to this study. The appendix contains the proof of
the technical lemmas
used in sections 4 and 6.  \\

\noindent {\bf Acknowledgements:} This paper was essentially
finished in the October of 2009 during a conference in honor of
Simon Donaldson at Northwestern University. With admiration, we
want to dedicate this modest paper to him for his teaching of
K\"ahler geometry to the first author in the last 12 years.
 Part of this work was done while both
 authors were visiting Stony Brook.  We wish to thank  both the
department of Mathematics and the Simons Center for Geometry and
Physics for their generous hospitality.  We also thank Professors
Blaine Lawson, Claude Lebrun, and Gang Tian for their interest in
this work. The second author would also like to thank Joel Fine,
Sean Paul and Zhan Wang for interesting discussions. Both authors
are partially supported by an NSF grant.

\section{The space of K\"ahler structures}
Here we review the infinite dimensional moment map picture
discovered by Fujiki(\cite{Fu}) and Donaldson(\cite{D1}). Let $(M,
\omega, J_0)$ be a compact K\"ahler manifold. Denote by $\J$ the
space of almost complex structures on $M$ which are compatible
with $\omega$, and by $\J^{int}$ the subspace of $\J$ consisting
of integrable almost complex structures compatible with $\omega$.
Then $\J$ is the space of smooth sections of an $Sp(2n)/U(n)$
bundle over $M$, so it carries a natural K\"ahler structure.
Indeed, there is a global holomorphic coordinate chart if we use
the ball model of the Siegel upper half space in the usual way.
$J_0$ determines a splitting $TM\otimes\C=T^{1,0}\oplus T^{0,1}$
such that $\omega$ induces a positive definite Hermitian inner
product on $T^{1,0}$, then $\J$ could be identified with the space
$$\OTS=\{\mu\in\OT|\A(\mu)=0, Id-\bar{\mu}\circ\mu>0\},$$ where
$\A$ is the composition
$\Omega^{0,p}(T^{1,0})\rightarrow\Omega^{0,p}(T^{*{0,1}})\rightarrow
\Omega^{0,p+1}$. An element $\mu$ corresponds to an almost complex
structure $J$ whose corresponding $(1,0)$ tangent space consists
of vectors of the form $X-\bar{\mu}(X)$($X\in T^{1,0}$).
$\J^{int}$ is a subvariety of $\J$ cut out by quadratic equations:
$$N(\mu)=\pb\mu+[\mu, \mu]=0.$$
 Denote by $\G$ the group of Hamiltonian diffeomorphisms of $(M,
\omega)$. Its Lie algebra is $C^{\infty}_0(M;\R)$. $\G$ will be
the infinite dimensional analogue of a compact group, though the
exponential map is not locally surjective for $\G$. $\G$ acts
naturally on $\J$, keeping $\J^{int}$ invariant. A.
Fujiki{\cite{Fu}} and S. Donaldson(\cite{D1}) independently
discovered that the $\G$ action has a moment map given by the
Hermitian scalar curvature functional
$S-\underline{S}$\footnote{Here $\underline{S}$ is the average of
scalar curvature, which indeed depends only on $[\omega]$ and
$c_1(\omega)$, not on the choice of any compatible $J$.}, which
can be viewed as an element in  $(C^{\infty}_0(M;\R))^*$ through
the $L^2$ inner product with respect to the measure $d\mu=\omega^n$.
When $J$ is integrable $S(J)$ is simply the Riemannian scalar
curvature of the Riemannian metric induced by $\omega$ and $J$. We
say $J_0\in \J$ is \emph{cscK} if $J_0$ is integrable and
$(\omega, J_0)$ has constant scalar curvature. So in the
symplectic theory we are naturally lead to  consider cscK metrics.

In the complex story, we need to look at $\G^{\C}$. Since $\G$ is
infinite dimensional, there may not exist a genuine
complexification $\G^{\C}$. Nevertheless, we can still define the
$\G^{\C}$ leaf of an integral complex structure $J_0$, as follows.
The infinitesimal action of $\G$ at a point $J\in \J$ is given by
$$\D_J: C^{\infty}_0(M;\R)\rightarrow \Omega^{0,1}_S(T^{1,0}); \phi\rightarrow \pb_{J}X_{\phi}.$$
This operator can be naturally complexified  to an operator from
$C^{\infty}_0(M;\C)=C^{\infty}_0(M;\R)\oplus
\sqrt{-1}C^{\infty}_0(M;\R)$ to $\OTS$. Then a complex structure
$J$ is on the $\G^{\C}$ leaf of $J_0$ if there is a smooth path
$J_t\in \J^{int}$ such that $\dot{J}_t$ lies in the image of
$\D_{J_t}$. $\G$ acts on the leaf naturally and the quotient is
the space of K\"ahler metrics cohomologous to $[\omega]_{J_0}$. So
the latter could be viewed as $``\G^{\C}/\G"$. We define the space
of K\"ahler potentials
$$\H=\{\phi\in C^{\infty}(M;\R)|\omega+\sqrt{-1}\p\pb\phi>0\}.$$
Then $\H/\R$ is formally the ``dual" symmetric space of $\G$. This
was made more precise by Mabuchi(\cite{M1}), Semmes(\cite{Se}) and
Donaldson(\cite{D2}). Define a Weil-Petersson type Riemannian
metric on $\H$ by
$$(\psi_1,\psi_2)_{\phi}=\int_M \psi_1\psi_2d\mu_{\phi}$$ for $\psi_1,\psi_2\in
T_{\phi}\H$. It can be shown that the Riemannian curvature tensor
is co-variantly constant and the sectional curvature is
non-positive. A path $\phi(t)$ in $\H$ is a \emph{geodesic} if it
satisfies the equation
$$\ddot{\phi}(t)-|\nabla_{\phi(t)}\dot{\phi}(t)|^2_{\phi(t)}=0.$$
 The first author(\cite{Ch1}) proved
the existence of a unique $C^{1,1}$ geodesic connecting any two
points in $\H$, and consequently that $\H$ is a metric space with
the distance given by the length of the $C^{1,1}$ geodesics. It is
proved in \cite{CC} that under this metric $\H$ is non-positively
curved in the sense of Alexanderov. So far the best regularity for
the Dirichlet problems of the geodesic equation was obtained by
Chen-Tian(\cite{CT}). The initial value problem for the geodesic
equation is in general not well-posed. But by the non-positiveness
of the curvature of $\H$, there should be lots of geodesic rays in
$\H$. In \cite{Ch3}, the first author proved the following general
theorem which we shall use later:

\begin{thm} Given a smooth geodesic ray $\phi(t)$ in $\H$ which is tamed by
 a bounded geometry, there is a unique relative $C^{1,1}$ geodesic
ray $\psi(t)$ emanating from any point $\psi$ in $\H$ such that
$$|\phi(t)-\psi(t)|_{C^{1,1}}\leq C.$$
\end{thm}

\begin{rmk} For the precise definition of ``tameness" we refer to \cite{Ch3}. But we point out that this is merely a technical condition imposed on the behavior of $\phi(t)$ at infinity so that the analysis on non-compact manifolds work. In our later applications where the geodesic ray $\phi(t)$ arises naturally from a test configuration with smooth total space, this assumption is always satisfied.
\end{rmk}

\begin{defi} Two geodesic rays $\phi(t)$ and $\psi(t)$ in $\H$ are
said to be \emph{parallel} if $$d_{\H}(\phi(t),\psi(t))\leq C.$$
\end{defi}
Hence it is clear by definition that if
$|\phi(t)-\psi(t)|_{C^0}\leq C$, then
$\phi$ and $\psi$ are parallel.\\

Analogous to the finite dimensional Kempf-Ness setting, there is a relevant
functional $E$ defined on $\H$, called the Mabuchi
\emph{K-energy}. It is the anti-derivative of the following closed
one-form:

\begin{equation}\label{Kenergy}
dE_{\phi}(\psi)=-\int_M(S(\phi)-\underline{S})\psi d\mu_{\phi}.
\end{equation} So the norm square of the gradient of $E$ is the Calabi energy:
$$Ca(\phi)=\int_M(S(\phi)-\underline{S})^2d\mu_{\phi}.$$ By a
direct calculation, along a smooth geodesic $\phi(t)$, we have
$$\frac{d^2}{dt^2}E(\phi(t))=\int_M|\D_t\dot{\phi}(t)|^2d\mu_{\phi(t)}\geq0.$$
According to \cite{Ch2}, $E$ can be extended to a continuous
function on all $C^{1,1}$ potentials in $\H$. However, it is not
clear why $E$ is still convex. The first author proved some weak
versions of convexity. In the case when $[\omega]$ is integral, we
gave simplified proofs in \cite{CS} using quantization(See also
\cite{Be}). We recall them for our later purpose.

\begin{lem}(\cite{Ch3}, \cite{CS}). \label{K energy convexity 1}
Given any $\phi_0$, $\phi_1\in \H$, we have
$$E(\phi_1)-E(\phi_0)\leq \sqrt{Ca(\phi_1)}\cdot d(\phi_0,\phi_1).$$
\end{lem}

\begin{lem}(\cite{Ch3}, \cite{CS})\label{K energy convexity} Given any  $\phi_0$, $\phi_1\in \H$, let $\phi(t)$ be
the $C^{1,1}$ geodesic connecting them. Then the derivatives of
$E(\phi(t))$ at the end-points are well-defined and they satisfy
the following inequality:
$$\frac{d}{dt}|_{t=0}E(
\phi(t))\leq \frac{d}{dt}|_{t=1}E(\phi(t)).$$
\end{lem}

This lemma implies that

\begin{lem}(\cite{CC})\label{Calabi decrease distance}
The Calabi flow on $\H$ decreases geodesic distance.
\end{lem}

\section{Lojasiewicz inequality}
In this section we recall Lojasiewicz's theory for the structure of a real analytic function. The following fundamental
structure theorem for real analytic functions is well-known:

\begin{thm}(Lojasiewicz inequality) \label{thmlo}Suppose $f$ is a real analytic function defined in a
neighborhood $U$ of the origin in $\R^n$. If $f(0)=0$ and $\nabla
f(0)=0$, then there exist constants $C>0$, and
$\alpha\in[\frac{1}{2}, 1)$, and shrinking $U$ if necessary,
depending on $n$ and $f$, such that for any $x\in V$, it holds
that
\begin{equation}\label{loja inequality}
|\nabla f(x)|\geq C\cdot|f(x)|^{\alpha}.
\end{equation}
\end{thm}

This type of inequality is crucial in controlling the behavior of
the gradient flow. If $\alpha=\frac{1}{2}$, then we get
exponential convergence. If $\alpha>\frac{1}{2}$, then we can
obtain polynomial convergence:

\begin{cor}\label{loja argu} Suppose $f$ is a non-negative real-analytic function defined in a
neighborhood $U$ of the origin in $\R^n$ with $f(0)=0$. Then there
exists a neighborhood $V\subset U$ of the origin such that for any
$x_0\in V$, the downward gradient flow of $f$:
\begin{equation*}
\left\{
 \begin{array}{ll}
          \frac{d}{dt}x(t)=-\nabla f(x(t)), \\
          x(0)=x_0.\\
\end{array}\right.
\end{equation*}
 converges uniformly to a limit $x_{\infty}\in U$ with
$f(x_{\infty})=0$. Moreover, we have the following estimate:

\begin{enumerate}

\item $$f(x(t))\leq C\cdot t^{-\frac{1}{2\alpha-1}};$$

\item $$d(x(t), x(\infty)) \leq C\cdot
t^{-\frac{1-\alpha}{2\alpha-1}},$$ where we assume the Lojasiewicz
exponent $\alpha> \frac{1}{2}$. 

\end{enumerate}
\end{cor}

\begin{proof} The proof is quite standard, and we call it ``Lojasiewicz arguments"
 for later reference. Denote $$V_{\delta}=\{x\in \R^n||x|\leq
\delta\},$$ and fix $\delta>0$ small so that inequality
$(\ref{loja inequality})$ holds for $x\in V_{\delta}$ . In our
calculation the constant $C$ may vary from line to line. If
$x(t)\in V_{\delta}$ for $t\in[0,T]$ , then we compute
$$\frac{d}{dt}f^{1-\alpha}(x(t))=-(1-\alpha)\cdot f^{-\alpha}(x(t))\cdot|\nabla f(x(t))|^2\leq -C\cdot|\dot{x}(t)|,$$
thus for any $T>0$,
$$\int_0^T|\dot{x}(t)|dt\leq \frac{1}{C}\cdot f^{1-\alpha}(x_0).$$
For any $\epsilon\leq \frac{\delta}{2}$ small, we choose
$\delta_2\leq \delta$ small such that $f(x)\leq
(C\cdot\epsilon)^{\frac{1}{1-\alpha}}$ for $x\in V_{\delta_2}$,
and $\delta_1=min\{\epsilon, \delta_2\}$, then the flow initiating
from any point $x_0\in V_{\delta_1}$ will stay in $V_{2\epsilon}$.
So the Lojasiewicz inequality holds for all $x(t)$. Now
$$\frac{d}{dt}f^{1-2\alpha}(x(t))=-(1-2\alpha)\cdot f^{-2\alpha}(x(t))\cdot |\nabla f(x(t))|^2\geq (2\alpha-1)\cdot C^2,$$
so $$f(x(t))\leq C\cdot t^{-\frac{1}{2\alpha-1}}.$$ For any
$T_1\leq T_2$, we get
$$d(x(T_1), x(T_2))\leq \int_{T_1}^{T_2}|\dot{x}(t)|dt\leq C \cdot {T_1}^{-\frac{1-\alpha}{2\alpha-1}}. $$
Therefore we obtain  polynomial convergence and the required
estimates.
\end{proof}

\section{Finite dimensional case}
\subsection{Kempf-Ness theorem}\label{Kempf-Ness}
Let $(M,\omega, J)$ be a K\"ahler manifold and  assume there is an
action of a compact connected group $G$ on $M$ which preserves the
K\"ahler structure. Let $\mu$ be the corresponding moment map.
This induces a holomorphic action of the complexified group
$G^{\C}$.  Then the Kempf-Ness theorem relates the complex
quotient by $G^{\C}$ to the symplectic
reduction by $G$(\cite{DK}).\\

\begin{thm}(Kempf-Ness)A $G^{\C}$-orbit contains a zero of the moment map if
and only if it is poly-stable. It is unique up to the action of
$G$. A $G^{\C}$-orbit is semi-stable if and only if its closure
contains a zero of the moment map; this zero is in the unique
poly-stable orbit in the closure of the original orbit.
\end{thm}

In this paper we are only interested in the uniqueness problem. We
will first give a proof in the finite dimensional case, using an
analytic approach. An essential ingredient in the proof of the
Kempf-Ness theorem is the existence of a function $E$, called the
\emph{Kempf-Ness function}. Given a point $x\in M$, one can define
a one-form $\alpha$ on $G^{\C}$ as:
$$\alpha_{g}(R_g\xi)=-\langle \mu(g.x), J\xi \rangle,$$
 where $R_g$ is the right translation by $g$ and $\xi\in \g_{\C}$. It is easy to check that $\alpha$ is closed and invariant under
 the left $G$-action. Then $\alpha$ is the pull back of a closed
 one-form $\bar{\alpha}$ from $G^{\C}/G$.
It is well known that $G^{\C}/G$ is always contractible, so $\alpha$
gives rise to a function $E$, up to an additive constant.
 Notice if the $G$ action is linearizable, this coincides
 with the usual definition given by the logarithm of the length of a vector on the induced line bundle.
 It is a standard fact that $E$ is
geodesically convex, i.e. $\bar{\alpha}$ is monotone along
geodesics in $G^{\C}/G$.  The critical points of $E$ consist
exactly of the zeroes of $\mu$ in the given $G^{\C}$ orbit. So any
$G^{\C}$ orbit contains at most one zero of the moment map, up to
the action of $G$. In the semi-stable case, we consider the
function $f(x)=|\mu(x)|^2$ on $M$, and its downward gradient flow
$x(t)$. The flow line is tangent to the $G^{\C}$ orbit and the
induced flow in $G^{\C}/G$ is exactly the downward gradient flow
of $E$.  We call either flow the \emph{Kempf-Ness flow}. As we
will see more explicitly later, a theorem of
Duistermaat(\cite{Le}) says that for $x(0)$ close to a zero of
$\mu$, the flow $x(t)$ converges polynomially fast to a limit in
$\mu^{-1}(0)$. Now suppose $x$ is semi-stable, and $x_1$, $x_2$
are two poly-stable points in $\overline{G^{\C}.x}$. W.L.O.G, we
can assume $\mu(x_1)=\mu(x_2)=0$. Take $y_1, y_2\in G^{\C}.x$ such
that $y_i$ is close to $x_i$. Then the gradient flows $x_i(t)$
converges to a point $z_i\in \mu^{-1}(0)$ near $x_i$. Denote by
$\gamma_i(t)$ the corresponding flow in $G^{\C}/G$. Since the
gradient flow of a geodesically convex function decreases the
geodesic distance, $d(\gamma_1(t), \gamma_2(t))$ is uniformly
bounded. By compactness, we conclude that $z_1$ and $z_2$ must be
in the same $G^{\C}$ orbit and by the uniqueness in the
poly-stable case, we see that $z_1$ and $z_2$ must lie in the same
$G$ orbit. By choosing $y_i$ arbitrarily close to $x_i$, we
conclude that $x_1$
and $x_2$ are in the same $G$ orbit.\\

The above argument proves the uniqueness of the poly-stable orbits
in the closure of a semi-stable orbit. There are technical
difficulties to extend this argument to the infinite dimensional
setting, due to the loss of compactness. As a result, we need to
investigate more about the gradient flow in the finite dimensional
case. What we shall show next
 is that for a semi-stable point, the gradient flow is
asymptotic to an ``optimal" geodesic ray at infinity.\\

\begin{figure}
 \begin{center}
  \psfrag{A}[c][c]{$\gamma(t)$}
  \psfrag{B}[c][c]{$\chi(t)$}
  \psfrag{C}[c][c]{$s$}
  \includegraphics[width=0.8 \columnwidth]{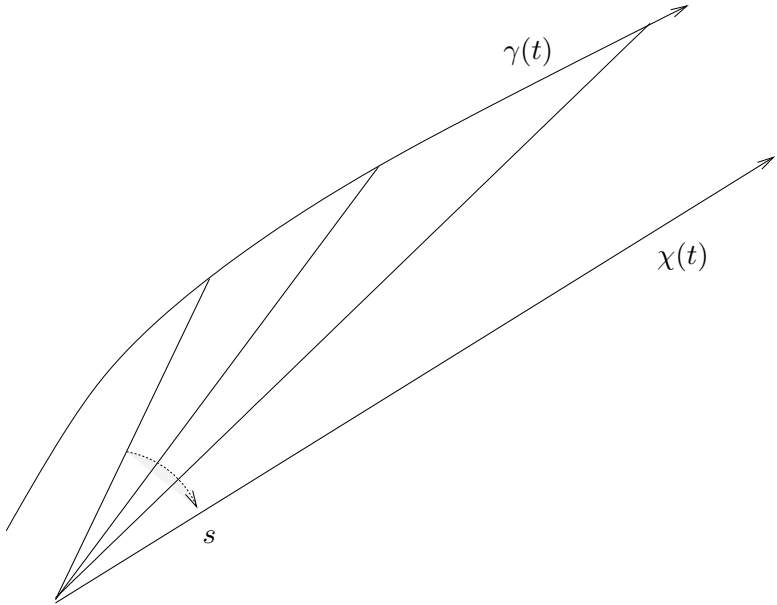}
  \caption{a curve asymptotic to a geodesic ray}
  \label{fig: parallel1}
 \end{center}
 \end{figure}

\begin{defi} \label{parallel curve}We say a curve $\gamma(t)(t\in [0,\infty))$ in a simply-connected
non-positively curved space is asymptotic to a geodesic ray
$\chi(t)$ if for any fixed $s>0$, $d(\gamma_t(s), \chi(s))$ tends
to zero as $t$ tends to $\infty$, where $\gamma_t$ is the geodesic
connecting $\chi(0)$ and $\gamma(t)$ which is parametrized by
arc-length. In other words, $\chi(t)$ is the point in the sphere
at infinity induced by $\gamma(t)$ as $t\rightarrow\infty$(see
figure \ref{fig: parallel1}).
\end{defi}

It follows from the definition that any  two geodesic rays
$\chi_1(t)$ and $\chi_2(t)$ that are both asymptotic to a given
curve $\gamma(t)$ must be \emph{parallel}, i.e. $d(\chi_1(t),
\chi_2(t))$ is
uniformly bounded.   \\

\subsection{Standard case}
Let $(V, J_0, g_0)$ be an $n$ dimensional unitary representation
of a compact connected Lie group $G$, so we have a group homomorphism:
$G\rightarrow U(n)$. $V$ is then a representation of the
complexified group $G^{\C}$. Denote by $\Omega_0$ the induced
K\"ahler form on $V$. It is easy to see that the $G$ action always
has a moment map $\mu: V\rightarrow \g^*\simeq \g$, where we have
identified $\g$ with $\g^*$ by fixing an invariant metric. It is
defined as
\begin{equation}
(\mu(v), \xi)=\frac{1}{2}\Omega_0(\xi.v, v).
\end{equation}
For any $v\in V$, denote the infinitesimal action of $G$ at $v$ by
$$L_v:\g\rightarrow V;\xi\mapsto \xi.v, $$
then it is easy to see that $$\mu(v)=\frac{1}{2}L_v^*(J_0v).$$
$L_v$ can also be viewed as a map from $\g^{\C}$ to $V$, and then
$\mu(v)=-\frac{1}{2}JL_v^*v$.\\

Now consider the function $f:V\rightarrow \R; v\mapsto
|\mu(v)|^2$, whose downward gradient flow equation is:
\begin{equation}\label{Flow}
\frac{d}{dt}v=-\nabla f(v)=-J_0L_v(\mu(v)).
\end{equation}

Since $f$ is a homogeneous polynomial, and thus real analytic, the
Lojasiewicz inequality holds for $f$, i.e. there exists constant
$C>0$ and $\alpha\in[\frac{1}{2},1)$, such that for $v$ close to
zero,
$$|\nabla f(v)|\geq C\cdot |f(v)|^{\alpha}.$$ The previous Lojasiewicz arguments show that for $v$ close to  $0$,
the flow (\ref{Flow}) starting from $v$ will converge
polynomially fast to a critical point of $f$.\\

From now on we assume $0$ de-stabilizes $v$, i.e. $0\in
\overline{G^{\C}.v} $. Thus the gradient flow (\ref{Flow})
converges to the origin by the uniqueness in the previous section.
Since everything is homogeneous, we can study the induced flow on
$\P(V)$. The action of $G$ is then holomorphic and Hamiltonian
with respect to the Fubini-Study metric on $\P(V)$, with moment
map $\hat{\mu}: \P(V)\rightarrow \g$. It is then easy to see that
 $$\hat{\mu}([v])=\frac{\mu(v)}{|v|^2}.$$
Let $\hat{f}=|\hat{\mu}|^2$, then we can study the downward
gradient flow of $\hat{f}$ on $\P(V)$:
\begin{equation}\label{Flowp}
\frac{d}{ds}[v]=-\nabla \hat{f}([v])=-J_0L_{[v]}(\hat{\mu}([v])).
\end{equation}
Let $\pi: V\rightarrow \P(V)$ be the quotient map, then clearly
$$\pi_*(\nabla f(v))=|v|^2\nabla \hat{f}([v]).$$

So the flow (\ref{Flowp}) is just a re-parametrization of the
image under $\pi$ of the flow (\ref{Flow}): if $v(t)$ satisfies
(\ref{Flow}), then $[v(s)]$ satisfies (\ref{Flowp}), with
$\frac{ds}{dt}=|v(t)|^2$. Since $\hat{f}$ is also real analytic,
the flow $[v(s)]$ converges polynomially fast to a unique limit
$[v]_{\infty}$.

\begin{lem}
$$\hat{\mu}([v]_{\infty})\neq 0$$.
\end{lem}
\begin{proof}
Otherwise $[v]$ is semi-stable with respect to the action of
$G^{\C}$ on $\P(V)$, thus the corresponding Kempf-Ness function
$\log |g.v|^2$ is bounded below on $G^{\C}$. This contradicts the
assumption that $0\in \overline{G^{\C}.v}$.
\end{proof}

 Thus we know that
$$\frac{\mu(v(s))}{|v(s)|^2}=\hat{\mu}([v]_{\infty})+O(s^{-\gamma})(\gamma>0)$$
is bounded away from zero when $s$ is large enough. So for $t$
sufficiently large, we have
$$|\nabla f(v(t))|^4\geq C\cdot|f(v(t))|^3.$$
The Lojasiewicz arguments then ensure that $v(t)$ actually
converges to $0$ in
the order $O(t^{-\frac{1}{2}})$. So we obtain $s\leq C\cdot \log t$.\\

Now since the gradient flow of $f$ is tangent to the $G^{\C}$
orbit, it can also be viewed as a flow on $G^{\C}/G$. This is
given by a path $\gamma(t)=[g(t)]$, where $g(t)\in G^{\C}$
satisfies
$$\dot{g}(t)g(t)^{-1}=-J\mu(g(t).v),$$ and the re-parameterized
path  corresponding to (\ref{Flowp}) is
$$\dot{g}(s)g(s)^{-1}=-J\hat{\mu}(g(s).[v]),$$
and
$$\frac{d}{ds}\gamma(s)=-J\hat{\mu}([v]_{\infty})+O(s^{-\gamma}).$$
In the following we shall use the re-parameterized version as
$|\frac{d}{ds}\gamma(s)|$  has a lower bound as
$s\rightarrow\infty$ which makes it more convenient to analyze the
asymptotic behavior.

\begin{thm}
$\gamma$ is asymptotic to a geodesic ray $\chi$ in $G^{\C}/G$.
Moreover, the direction of $\gamma$ is conjugate to
$\frac{\hat{\mu}([v]_{\infty})}{|\hat{\mu}([v]_{\infty})|}$ under
the adjoint action of $G$.
\end{thm}
\begin{proof} We already know $\dot{\gamma}(s)$ is getting close
to $\hat{\mu}([v]_{\infty})$, but this is not sufficient to
conclude that $\gamma$ is asymptotic to a geodesic ray with
direction $\hat{\mu}([v]_{\infty})$. We shall analyze this more carefully, by elementary geometry. First it is easy to see that
$$|\ddot{\gamma}(s)|=|L_{[v](s)}^*L_{[v](s)}\hat{\mu}([v](s))|,$$
where $L_{[v](s)}$ is the infinitesimal action of $\g$ at
$[v](s)$. Since $[v](s)\rightarrow [v]_{\infty}$ as $s\rightarrow
\infty$, by corollary \ref{loja argu} we get
$$\int_{t}^{\infty}|\ddot{\gamma}(s)|ds\leq C\int_t^{\infty}|L_{[v](s)}\hat{\mu}([v](s))|ds=C \int_t^{\infty}
|\nabla \hat{f}(s)|ds\leq C\cdot t^{-\beta},
$$
where $\beta=\frac{1-\alpha}{2\alpha-1}>0$. Notice that here $\alpha$ is the exponent appearing in the Lojasiewicz inequality for $\hat{f}$, not the original $f$. From the above we know
$\lim_{s\rightarrow\infty}|\dot{\gamma}(s)|=
|\hat{\mu}([v]_{\infty})|>0$, so if we parameterize $\gamma$ by
arc-length and denote the resulting path by  $\tgamma(u)$, then we
have
$$|\ddot{\tgamma}(u)|=|\dot{\gamma}(s)|^{-2}|\ddot{\gamma}(s)
-\frac{\langle\ddot{\gamma}(s),
\dot{\gamma}(s)\rangle}{|\dot{\gamma}(s)|^2}\dot{\gamma}(s) |\leq
C\cdot |\ddot{\gamma}(s)|.
$$ Therefore 
$$\int_{t}^{\infty}|\ddot{\tgamma}(u)|du\leq C\cdot t^{-\beta},
$$
Now for any $u>0$, let $\tgamma_u(v)$($v\in[0,1$) be the geodesic in $G^{\C}/G$
connecting $\tgamma(0)$ and $\tgamma(u)$. Denote by $L_u(v)$($v\in
[0, u]$) the distance between $\tgamma(v)$ and $\tgamma_u(v)$.
Then $L_u(0)=L_u(u)=0$ and a standard calculation of the second
variation of length(using the non-positivity of the sectional curvature of
$G^{\C}/G$) gives
$$\frac{d^2}{dv^2}L_u(v)\geq -|\ddot{\tgamma}(v)|.$$
Now define the function
$$f_u(v)=\int_0^v\int_{w}^{\infty}|\ddot{\tgamma}(r)|drdw-\frac{v}{u}\int_0^u\int_w^{\infty}|\ddot{\tgamma}(r)|drdw.$$
Then it is well-defined by the decay of $|\ddot{\tgamma}|$, and
$f_u(0)=f_u(u)=0$ and
$$\frac{d^2}{dv^2}f_u(v)=-|\ddot{\tgamma}(v)|.$$
Thus by maximum principle $L_u(v)\leq f_u(v)$ for all $u>0$ and
$v\in[0,u]$. Fix $v$ we see
$$\sup_{u}L_u(v)\leq\int_0^v\int_{w}^{\infty}|\ddot{\tgamma}(r)|drdw\leq C\cdot v^{1-\beta}.$$
Moreover, for any $u_2>u_1>>1$, by  comparison argument the angle
between $\tgamma_{u_1}$ and $\tgamma_{u_2}$ is bounded by
$d(\tgamma_{u_1}(u_1), \tgamma_{u_2}(u_1))/u_1=L_{u_2}(u_1)/u_1$,
which is controlled by $C\cdot u_1^{\beta-1}$. Thus we conclude
that the direction of $\tgamma_u$ is converging uniformly to some
limit direction and so $\tgamma$(and thus $\gamma$) is asymptotic
to a geodesic ray $\chi$ starting from $\gamma(0)$. Now for any
$s>0$ by the same way we get a geodesic ray $\chi_s$ starting from
$\gamma(s)$ which is asymptotic to $\gamma$. So the rays $\chi_s$
are all asymptotic to each other and one could easily see that
they are all parallel, and then $\dot{\chi}_s(0)$ are all
conjugate to each other under the action of $G$. On the other
hand, if we denote by $\gamma_{s,t}(u)$($u\in[0,1]$) the geodesic
connecting $\gamma(s)$ and $\gamma(t)$ for $s<t$, then again by
second variation,
$$\frac{d}{dt}\langle\frac{\dot{\gamma}(t)}{|\dot{\gamma}(t)|}, \frac{\dot{\gamma}_{s,t}(1)}{|\dot{\gamma}_{s,t}(1)|}\rangle\geq -C\frac{|\ddot{\gamma}(t)|}{|\dot{\gamma}(t)|}\geq -C|\ddot{\gamma}(t)|.$$
So we get
$$\langle\frac{\dot{\gamma}(t)}{|\dot{\gamma}(t)|}, \frac{\dot{\gamma}_{s,t}(1)}{|\dot{\gamma}_{s,t}(1)|}\rangle\geq 1-\int_s^t |\ddot{\gamma}(u)|du\geq 1-C\cdot s^{-\beta}.$$
We know $\dot{\gamma}(t)=J\hat{\mu}([v]_{\infty})+O(t^{-\alpha})$,
and as $t\rightarrow\infty$ up to the adjoint action of $G$  we have
$$\frac{\dot{\gamma}_{s,t}(1)}{|\dot{\gamma}_{s,t}(1)|}\rightarrow\dot{\chi}(s).$$
So let $s\rightarrow\infty$ we see $\dot{\chi}(0)$ is conjugate to
$\frac{\hat{\mu}([v]_{\infty})}{|\hat{\mu}([v]_{\infty})|}$ under
the adjoint action of $G$.
\end{proof}

From the proof of the above theorem we see that $\chi(s)$ also
degenerates $v$ to the origin since the path $v(t)$ is of order
$O(t^{-\frac{1}{2}})=O(e^{-C\cdot s})$. By Kempf(\cite{Ke}) and
Ness(\cite{Ne}), the direction $\hat{\mu}([v]_{\infty})$
is indeed rational, i.e. it generates an algebraic one-parameter subgroup
$\lambda:\C^*\rightarrow G^{\C}$. Moreover, the direction
$\hat{\mu}([v]_{\infty})$ is the unique(up to
the adjoint action of $G$) optimal direction for $v$ in the sense of
Kirwan(\cite{Ki})) and
Ness(\cite{Ne}).\\

\subsection{Linear Case} Now we suppose $G$ acts linearly on $(V=\C^n,\Omega,J_0)$
 where $J_0$ is the standard complex structure on $\C^n$ and $\Omega$ is a real-analytic symplectic form compatible
 with $J_0$. Then the action has a
real-analytic moment map $\mu$ with $\mu(0)=0$. $\mu$ is not necessarily  standard
but the
 Lojasiewicz inequality still holds for $f=|\mu|^2$. Suppose $0\in \overline{G^{\C}.v}$, then the downward gradient flow
 $v(t)$ of $f(v)=|\mu(v)|^2$ converges to the origin polynomially
 fast. Let $\hat{v}(t)$ be the downward gradient flow of $\hat{f}(v)=|\hat{\mu}(v)|^2$, where
 $\hat{\mu}$ is the moment map for the linearized $G$ action on $(V=T_0V,\Omega_0,
 J_0)$. By the arguments in the previous section, $\hat{v}(t)$
 converges to zero in the order $O(t^{-\frac{1}{2}})$ and the
 corresponding flow $\hat{\gamma}(t)$ is asymptotic to a rational
 geodesic ray $\chi(t)$. Let $\gamma(t)$ in
 $G^{\C}/G$ be the flow corresponding to $v(t)$, we want to
 show $\gamma(t)$ is also asymptotic to $\chi(t)$. It suffices to
 bound the distance $L(t)$ between $\gamma(t)$ and
 $\hat{\gamma}(t)$. Let $\psi_t(s)(s\in[0,1])$ be the geodesic
 connecting $\gamma(t)$ and $\hat{\gamma}(t)$, then
 \begin{eqnarray*}
 \frac{d}{dt}L(t)&=& \frac{1}{L(t)}\langle\dot{\psi}(1), \hat{\mu}(\hat{v}(t))\rangle-\frac{1}{L(t)}\langle\dot{\psi}(0),
 \mu(v(t))\rangle \\&=& \frac{1}{L(t)}(\langle\dot{\psi}(1),
 \mu(\hat{v}(t))\rangle-\langle\dot{\psi}(0),
 \mu(v(t))\rangle)+\frac{1}{L(t)}\langle\dot{\psi}(1),
 \hat{\mu}(\hat{v}(t))-
 \mu(\hat{v}(t))\rangle\\&\leq& |\hat{\mu}(\hat{v}(t))-
 \mu(\hat{v}(t))|,
 \end{eqnarray*}
 where we used the fact that the Kempf-Ness function is
 geodesically convex.
To estimate the last term, notice since the $G$ action is linear,
we have for any $\xi\in\g$ 
\begin{eqnarray*}
\langle\mu(v), \xi\rangle&=&\langle\mu(0)+\int_0^1\frac{d}{dt}\mu(tv)dt, \xi\rangle\\
&=&\int_0^1 \Omega_{tv}({\xi}.tv, v)dt\\
&=&\frac{1}{2}\Omega_0({\xi}.v, v)dt+O(|v|^3)\\
&=&\langle\hat{\mu}(v), \xi\rangle+O(|v|^3).
\end{eqnarray*}
 From the previous secion
we know $\hat{v}(t)=O(t^{-\frac{1}{2}})$, so we obtain
$$\frac{d}{dt}L(t)\leq C\cdot t^{-\frac{3}{2}},$$ and so
$L(t)$ is uniformly bounded. Therefore, we conclude the following
theorem:
\begin{thm}
Suppose $G$ acts Hamiltonian linearly  on $(V, \Omega, J_0)$, with
the moment map given by $\mu$. Suppose also a vector $v_0$ is
de-stabilized by the origin. Let $v(t)$ be the downward gradient
flow of $|\mu|^2$ emanating from $v$, then $v(t)$ converges to $0$
in the order $O({t^{-\frac{1}{2}}})$. Let $\gamma(t)$ be the
corresponding flow in $G^{\C}/G$, then there exists a geodesic ray
$\chi$ in $G^{\C}/G$, which is asymptotic to  $\gamma$. Moreover,
$\chi$ is rational.
\end{thm}

\subsection{General Case} In general we need to linearize the
problem, using the Marle-Guillemin-Sternberg normal form. Let $(M,
\omega, J, G, \mu)$ be a real analytic Hamiltonian $G$-action on a
real analytic K\"ahler manifold. Choosing a bi-invariant metric on
$\g$ we can identify $\g$ with $\g^*$. Suppose $x\in M$ is a zero
of $\mu$. Let $G_0$ be the isotropy group of $x$ and $\g_0$ be its
Lie algebra.  The bi-invariant product on $\g$ allows a $G_0$
invariant splitting:
$$\g=\g_0\oplus \m.$$
Notice $\g.x\subset (\g. x)^{\omega}$. Denote by $N$ the
orthogonal complement of $\g.x\oplus J\g.x$ in $T_xM$, then $N$ is
$G_0$-invariant and the linear $G_0$ action on $N$ has a canonical
moment map $\mu_N: N\rightarrow \g_0$. Let
$$Y=G\times_{G_0}(\m\times N),$$ then $G$ acts naturally on $Y$
on the left.

\begin{lem}(Marle-Guillemin-Sternberg \cite{GS}, \cite{OR})
\label{MGS} There exists a symplectic form $\omega$ defined in a
neighborhood $U$ of $[e,0,0]$ in $Y$, under which the $G$ action
is Hamiltonian with a moment map given  by
$$\mu: U\rightarrow \g; [g, \rho, v]\rightarrow Ad_{g}^*(\mu_N(v)+\rho).$$
There exists a local $G$ equivariant symplectic diffeomorphism
$\Phi: Y\rightarrow M$ which respects the moment maps, and
satisfies $\Phi([e,0,0])=x$, $\Phi^*J-J_0=O(r^2))$ on $N$ and
$\Phi^*J=J_0$ at $[e,0,0]$. Here $J_0$ is the canonical
$G$-invariant almost complex structure on $Y$ induced by $J$,
which will be more explicit in the proof. Moreover, we can take
$\Phi$ to be real analytic if everything we start with is so.
\end{lem}
The only new feature here is the control on the complex structure.
The proof of this theorem is a bit technical and will be deferred
to the appendix.\\

 From now on we will work on $(U, \Omega_0, J)$ where we also denote by $J$ the pullback $\Phi^*J$.

\begin{thm}\label{thmfinite}
Suppose $y\in U$ is de-stabilized by $x$, then the Kempf-Ness flow
$y(t)$ of $|\mu|^2$ converges to $y_{\infty}\in G.x$ polynomially
fast. Moreover the corresponding  flow $\gamma(t)$ in $G^{\C}/G$
is asymptotic to a geodesic ray $\chi(t)$ which is rational and
also degenerates $y$ to $y_{\infty}$.
\end{thm}

\begin{rmk}
Here we could define $\chi(t)$ as the ``optimal" degeneration of
$y$, generalizing the usual definition in the linear case.
\end{rmk}

To prove the theorem, we study the function $f=|\mu|^2$ on $U$. By
definition,
$$f([g, \rho, v])=|\rho|^2+|\mu_{N}(v)|^2,$$
$$\nabla f([g,\rho, v])=J[L_g\rho, ad_{\mu_N(v)}\rho, \mu_N(v).v],$$
Since $f$ is real analytic, we have for some $\alpha\in
[\frac{1}{2},1)$ that
$$|\nabla f|\geq C\cdot|f|^{\alpha}.$$
Therefore $y(t)$ converges to a zero $y_{\infty}$ of $\mu$
polynomially fast. By uniqueness, $y_{\infty}\in G.x$. Without
loss of generality, we will assume $y_{\infty}=x$ from now on, and
we shall distinguish between two cases.\\

In the first case we assume $G_0=G$, then $\m=0$, and we are
essentially reduced to the linear case. What we obtain is a
K\"ahler manifold $(U\subset N,\Omega_0, J)$. We just need to
holomorphically
linearize the $G$ action:\\

\begin{lem}\label{average trick}
There exits a $G$-equivariant holomorphic embedding
$$\Phi: (V\subset T_0U, J_0)\hookrightarrow (U, J); 0\mapsto x.$$
\end{lem}

\begin{proof} Shrinking $U$ if necessaray, we can first choose a holomorphic embedding
$$\Psi: (U, J)\hookrightarrow (T_0U, J_0); x\mapsto 0 .$$ Again Shrinking $U$ if necessary, define
$$\hat{\Psi}: (U, J)\rightarrow (T_0U, J_0); y\mapsto \frac{1}{|G|}\int_G g^{-1}\cdot \Psi(g.y)d\mu,$$
where $\mu$ is a Harr measure on $G$. Then $\hat{\Psi}$ is
holomorphic, and $d\hat{\Psi}_x=d\Psi_x$, so $\hat{\Psi}$ is an
embedding near $x$. Then we can just take $\Phi=\hat{\Psi}^{-1}$.
\end{proof}

Now using $\Phi$ we can work on $(V_1, \Omega=\Phi^*\Omega_0,
J_0)$ with a linear Hamiltonian of $G$, and the linear theory
in the previous section applies to conclude the theorem in this case.\\

In the second case we assume $G_0$ is a proper subgroup of $G$. We
will try to reduce to the first case. It is easy to see that the
$G_0$ action on $Y$ is also Hamiltonian, with a moment map
$\hat{\mu}$ equal to the orthogonal projection of $\mu$ to $\g_x$.
Therefore,
$$\hat{\mu}([g, \rho, v])=Ad^*_g\mu_N(v).$$

Denote by $G_0^{\C}$ the isotropy group of $x$.
\begin{lem}
$G_0^{\C}$ is the complexification of $G_0$(hence is reductive).
\end{lem}
\begin{proof}This lemma is well-known. In the Lie algebra level,
we just need to show if $\xi.x+J\eta.x=0$  for some $\xi, \eta\in
\g$, then $\xi.x=\eta.x=0$. This follows easily from the
definition of the moment map:
$$\omega(\eta.x, J\eta.x)=(d\mu(J\eta.x), \eta)=(d\mu(J\eta.x+\xi.x), \eta)-(Ad_{\xi}^*\mu(x),\eta)=0.$$
Hence $\eta.x=0$ and $\xi.x=0$.
\end{proof}

\begin{lem} We can  choose a point in the $G^{\C}$ orbit of $y$,
 denoted by $\hat{y}$, so that $x$ de-stabilizes $\hat{y}$ for the
group $G_0$. \end{lem}

\begin{proof}

It suffices to find $\hat{y}$ in the $G^{\C}$ orbit of $y$ such
that $x$ lies in the closure of $G_0^{\C}.\hat{y}$. To do this, we
first choose an arbitrary holomorphic map $\Psi: T_xM\rightarrow
M$ with $\Psi(0)=x$ and $d\Psi(0)=Id$. As before we can linearize
the action so that $\Psi$ is $G_0$-equivariant. $T_xM$ has a
$\C$-linear decomposition
$$T_xM=\g^{\C}.x\oplus N,$$
where $N$ is as before the orthogonal complement of
$\g^{\C}.x=\g.x\otimes\C=\g.x\oplus J_0(\g.x)$.
  Then we define
$$\Phi:G^{\C}\times_{G^{\C}_0}N\rightarrow M; [(g, v)]\rightarrow
g. \Psi(v).$$ This is a local diffeomorphism around $[(Id,0)]$. So
for any $y$ close to $x$, there is a unique $(g, v)\in
G^{\C}\times N$ which is close to $[(Id,0)]$  such that
$y=g.\Psi(Id,v)$. Let $\hat{y}=\Psi(Id, v)$. We claim $x\in
\overline{G_0^{\C}.\hat{y}}$. Notice that the Kempf-Ness flow
$y(t)$ converges to $x$, so this gives rise to a smooth family
$(g(t), v(t))$ with $y(t)=g(t).\Psi(Id, v(t))$. Let
$\hat{y}(t)=\Psi(Id, v(t))$. Since $y(t)$ all lie in the same
$G^{\C}$ orbit, so are $\hat{y}(t)$. Thus all $v(t)$ lie in the
$G^{\C}_0$ orbit of $v$, and
$$\lim_{t\rightarrow\infty}v(t)=0.$$ Therefore, $x\in
\overline{G_0.\hat{y}}$.
\end{proof}

 Let $\tilde{y}(t)$ be the downward gradient flow of $f$ with $\tilde{y}(0)=\hat{y}$,
 and $\hat{y}(t)$ be the downward gradient flow of
$\hat{f}=|\hat{\mu}|^2$ with $\hat{y}(0)=\hat{y}$. Let
$\tilde{\gamma}(t)$ and $\hat{\gamma}(t)$ be the corresponding
path in $G^{\C}/G$ and $G_0^{\C}/G_0$ respectively. Then  the
previous linear theory tells that $\hat{y}(t)$ converges to $x$ in
the order $O(t^{-\frac{1}{2}})$ and $\hat{\gamma}(t)$ is
asymptotic to a rational geodesic ray $\chi(t)$ with the same
degeneration limit. On the other hand $G_0^{\C}/G_0$ is naturally
a totally geodesic submanifold of $G^{\C}/G$, and next we
will prove that the distance between $\tilde{\gamma}(t)$ and $\hat{\gamma}(t)$ in $G^{\C}/G$ is uniformly bounded.\\

We denote by $\psi_t(s)(s\in[0,1])$ the geodesic in $G^{\C}/G$
connecting $\tilde{\gamma}(t)$ and $\hat{\gamma}(t)$, and $L(t)$
the length of $\psi_t$, then it is easy to see that
\begin{eqnarray*}
\frac{d}{dt}L(t)&=&\frac{1}{L(t)}(\mu(y(t)),\dot{\psi}_t(0))-\frac{1}{L(t)}(\hat{\mu}(\hat{y}(t)),\dot{\psi}_t(1))\\\\&=&
\frac{1}{L(t)}(\mu(y(t)),\dot{\psi}_t(0))-\frac{1}{L(t)}(\mu(\hat{y}(t)),\dot{\psi}_t(1))+
\frac{1}{L(t)}(\mu(\hat{y}(t))-\hat{\mu}(\hat{y}(t)),\dot{\psi}_t(1))\\\\
&\leq& |\mu(\hat{y}(t))-\hat{\mu}(\hat{y}(t))|,\\
\end{eqnarray*}
where again we have used the convexity of the Kempf-Ness function.
In our situation, $\mu-\hat{\mu}=Ad_g^*\rho$. Here $g(t)$ are
$\rho(t)$ are uniquely determined by the choice at $t=0$ if we
require $\dot{\rho}(t)\in\m$ and $g(t)^{-1}\dot{g}(t)\in \m$. Now
at $\hat{y}(t)=[g(t), \rho(t),v(t)]$, we have
\begin{eqnarray*}
\nabla \hat{f} &=&
J.([0,ad_{\mu_N(v)}^*\rho,\mu_N(v).v])\\\\&=&[ad_{\mu_N(v)}^*\rho,
0,
J_0\cdot(\mu_N(v).v)]+(J-J_0)ad_{\mu_N(v)}^*\rho+(J-J_0)\mu_N(v).v.
\end{eqnarray*}
Therefore,
\begin{eqnarray*}|\frac{d}{dt}\rho(\hat{y}(t))|&=&|\Pi_{\m}(\nabla
\hat{f})|\\&\leq& C\cdot|J-J_0||\mu_N(v)||\rho|+C\cdot
d(\hat{y}(t), x)^2|\mu_N(v).v|)\\&\leq&
C\cdot(t^{-\frac{3}{2}}|\rho|+t^{-\frac{5}{2}}).
\end{eqnarray*} Since $\rho(\infty)=0$, we first get
$$|\rho(t)|\leq C \cdot t^{-\frac{1}{2}}.$$ Then plug back into the
previous inequality and repeat to obtain
$$\frac{d}{dt}\rho(\hat{x}(t))\leq C\cdot
t^{-\frac{5}{2}},$$ and then
$$|\rho(\hat{x}(t))|\leq C\cdot t^{-\frac{3}{2}}.$$
 So $$ L(t)\leq \int_{1}^t
s^{-\frac{3}{2}}ds+C\leq C.$$ Therefore $L(t)$ is uniformly
bounded.

By definition, we see that $\tilde{\gamma}(t)$ is also asymptotic
to the geodesic ray $\chi(t)$. Now the original $\gamma(t)$ is
also asymptotic to $\chi(t)$ again because that the Kempf-Ness
flow in $G^{\C}/G$ decreases the geodesic distance.

Then it is easy to see that $\chi(t)$ has the same degeneration
limit as $\gamma(t)$. So this completes the proof of theorem
\ref{thmfinite}.

\section{Stability of the Calabi flow}
 We first recall the definition of the Calabi flow. It is an
infinite dimensional analogue of the previously mentioned
Kempf-Ness flow. Let $(M, \omega, J_0)$ be a K\"ahler manifold. As
before, we have the group $\G$ acting on $\J$ and preserves
$\J^{int}$. The action of $\G$ on $\J$ has a moment map given by
the Hermitian scalar curvature functional
$$S-\underline{S}: \J\rightarrow C^{\infty}_0(M;\R).$$ Its norm is called the \emph{Calabi functional}:
 $$Ca(J)=\int_M (S(J)-\underline{S})^2d\mu_{\omega}.$$ The gradient of $Ca$ under the natural metric on $\J$
 is given by
 \begin{equation*}
 \nabla Ca(J)=\frac{1}{2}J\D_JS(J).\footnote{The factor comes from the fact that the metric we choose on $\J$
is $(\mu_1, \mu_2)_J:=2Re\int_M\langle\mu_1,
\mu_2\rangle_J\omega^n$.}
 \end{equation*} The \emph{Calabi flow} is the downward gradient flow of $Ca$
 on $\J^{int}$.
 Its equation is given by
\begin{equation} \label{Calabi flow J}
\frac{d}{dt}{J(t)}=-\frac{1}{2}J(t)\D_{J(t)}S(J(t)).
\end{equation} As in the finite dimensional space, the
Calabi flow can be lifted to $\G^{\C}/\G$, which in this case is
just the space of K\"ahler metrics
$$\H_{J}=\{\phi\in
C^{\infty}_0(M;\R)|\omega+\sqrt{-1}\p_J\bar{\p}_J\phi>0\}.$$ The
equation reads:
\begin{equation} \label{Calabi flow H}
\frac{d}{dt}\phi(t)=S(\phi(t))-\underline{S}.
\end{equation}
By (\ref{Kenergy}), this is also the downward gradient flow of the
Mabuchi functional $E$. The two equations (\ref{Calabi flow J})
and (\ref{Calabi flow H}) are essentially equivalent:

\begin{lem}\label{Two Calabi flows}
Any solution of (\ref{Calabi flow H}) naturally gives rise to a
solution of (\ref{Calabi flow J}); any solution $J(t)$ of
(\ref{Calabi flow J}) induces a solution of (\ref{Calabi flow H}),
if $J(t)$ all lie in $\J^{int}$.
\end{lem}
\begin{proof} Given a path $\phi(t)\in \H$, we consider the time-dependent vector fields
$X(t)=-\frac{1}{2}\nabla_{\phi(t)}\dot{\phi}(t)$. Let $f_t$ be the
family of diffeomorphisms generated by $X(t)$. Then
$f_t^*(\omega+\sqrt{-1}\p\pb\phi(t))=\omega$. Let $J(t)=f_t^*J$.
Then
$$\frac{d}{dt}J(t)=-\frac{1}{2}J(t)\D_{J(t)}\dot{\phi}(t).$$ This proves the first half of the lemma. For the
second half, if $J(t)$ is a solution to (\ref{Calabi flow J}). We
again consider the vector fields
$X(t)=\frac{1}{2}\nabla_{J(t)}S(J(t))$ and the induced
diffeomorphisms $f_t$. Then $f_t^*J(t)=J(0)$ since $J(t)\in
\J^{int}$, and
$f_t^*\omega=\omega+\sqrt{-1}dJ(0)d\phi(t)$, with
$\frac{d}{dt}\phi(t)=S(\phi(t))-\underline{S}$.
\end{proof}

Equation (\ref{Calabi flow J}) is not parabolic, due to the $\G$
invariance.  But (\ref{Calabi flow H}) is parabolic and we have
the following estimates:
\begin{lem}\label{lem2}(see \cite{CH2}) Suppose there are constants $C_1, C_2>0$ such that along the Calabi flow:
\begin{equation}
\left\{
 \begin{array}{ll}
\frac{\p \phi}{\p t}=S-\underline{S} \\
          \phi(0)=\phi_0,\\
\end{array}\right.
\end{equation}
 we have $$||Rm(g(t))||_{L^{\infty}(g(t))}\leq C_1,$$ and the Sobolev constant of $g(t)$ is bounded by $C_2$ for all $t\in[0, T)$, then for
any $l>0$, and $t\in [1, T)$, we have $$||\nabla_t^l
Rm(g(t))||_{L^{\infty}(g(t))}\leq C, $$ where $C>0$ depends only
$C_1, C_2, l, n$.
\end{lem}
The Calabi flow equation in the form (\ref{Calabi flow H}) was
first proposed by E. Calabi(\cite{Ca1}, \cite{Ca2}) to find
extremal metrics in a fixed K\"ahler class. The short time
existence was established
by Chen-He(\cite{CH1}). They also proved the global existence assuming Ricci curvature bound. \\

The equation (\ref{Calabi flow J}) also has its own advantage.
Namely,  when the space $\H$ does not admit any cscK metric, the
solution of  equation $(\ref{Calabi flow H})$ must diverge when
$t\rightarrow\infty$. However, it is still possible that the
corresponding $J(t)$ still converges in the bigger ambient space
$\J$. In this section we are interested in the Calabi flow
(\ref{Calabi flow J}) starting from an integrable complex
structure in a neighborhood of  a cscK metric. We shall prove the
following theorem:\\

\begin{thm} \label{stabilityinfi}Suppose $J_0\in \J$ is  cscK. Then
there exists a small $C^{k,\lambda}(k\gg1)$ neighborhood $\U$ of
$J_0$ in $\J^{int}$, such that the Calabi flow $J(t)$ starting
from any $J\in \U$ exists globally and converges  polynomially
fast to a cscK metric $J_{\infty}\in \J$ in $C^{k,\lambda}$
topology. Up to a Hamiltonian diffeomorphism we can assume
$J_{\infty}$ is smooth, then the convergence is also in
$C^{\infty}$.
\end{thm}

\begin{rmk} When $J$ lies on the leaf of $J_0$, i.e. the
corresponding K\"ahler metrics are in the same K\"ahler classes,
this was proved in \cite{CH1} and the convergence is indeed
exponential. In general, the convergence is exponential if and
only if $J_0$ and $J_{\infty}$ are on the same $\G^{\C}$ leaf.
\end{rmk}

\begin{rmk} There are also studies of stabiliy of other geometrical flows (such as K\"ahler-Ricci flow) in
K\"ahler geometry when the complex structure is deformed, see for
example \cite{CLW}, \cite{TZ2}... We believe the idea in this
section could also apply to other settings. In a sequel to this
paper(\cite{SW}), the second author and Y-Q. Wang proved a similar
stability theorem for the K\"ahler-Ricci flow on Fano manifolds.
We should mention that two alternative approaches in the study of
the stability of K\"ahler-Ricci flow have been announced by
C.Arezzo-G. La Nave and G. Tian-X. Zhu.
\end{rmk}

In general this type of stability result is based on a very rough
a priori estimate of the length of the flow and the parabolicity.
Here the key ingredient is the following Lojasiewicz type
inequality which yields the required a priori estimate.\\

\begin{thm}\label{thm1}Suppose $J_0\in\J^{int}$ is cscK, then there exists
a $L^2_k$($k\gg1$) neighborhood $\U$ of $J_0$ in $\J^{int}$ and
constants $C>0$, $\alpha\in [\frac{1}{2},1)$ such that for any
$J\in \U$, the following inequality holds:

\begin{equation}\label{lojainfi}
||\D_J S(J)||_{L^2}\geq C \cdot
||S(J)-\underline{S}||^{2\alpha}_{L^2}, \end{equation} where $\D_J
\phi=\pb_{J} X_{\phi}+\bar{X}_{\phi}.N_J$.
 When $J$ is integrable, $\D_J\phi=\pb_J X_{\phi}$ is the
Lichnerowicz
operator.\\
\end{thm}

\begin{rmk} The Lojasiewicz inequality was first used by L.
Simon(\cite{Si}) in the study of convergence of parabolic P.D.E's.
R\aa de(\cite{Ra}) used Simon's idea to study the convergence of
the Yang-Mills flow on two or three dimensional manifold. It also
appeared in the study of asymptotic behavior in Floer theory in
\cite{D4}. Here we
follow \cite{Ra} closely.\\
\end{rmk}

  We begin the proof by reducing the problem to a
finite dimensional
one and then use Lojasiewicz's inequality(theorem \ref{thmlo}).\\

To simplify the notation, we assume the function spaces appearing
below consist of normalized functions, i.e. functions with average
zero. We have the elliptic complex at $J_0$(see \cite{FS}):
$$L^2_{k+2}(M;
\C)\stackrel{\D_0}{\longrightarrow}T_{J_0}\J=L^2_k(\Omega_S^{0,1}(T^{1,0}))\stackrel{\pb_0}{\longrightarrow}
L^2_{k-1}(\Omega_S^{0,2}(T^{1,0})), $$ where
$\Omega_S^{0,p}(T^{1,0})$ is the kernel of the operator $\A$ in
section 3.  So we have an $L^2$ orthogonal decomposition:
$$\Omega_S^{0,1}(T^{1,0})=Im \D_0\oplus Ker \D_0^*. $$ On the other hand, the infinitesimal action of the gauge group
$\G$ is just the restriction of $\D_0$ to $L^2_{k+2}(M;\R)$, which
we denote by $Q_0$. Since $J_0$ is cscK, $\D_0^*\D_0$ is a real
operator. Thus
$$Im(\D_0)=\D_0(L^2_{k+2}(M; \R))\oplus \D_0(L^2_{k+2}(M;
\sqrt{-1}\R))$$ is an $L^2$ orthogonal decomposition, so
$$L^2_k(\Omega_S^{0,1}(T^{1,0}))=Im Q_0\oplus Ker Q_0^*,$$
where explicitly, $Q_0^*\mu=Re \D_0^*\mu$.\\

Now as in section 2 we identify a $L^2_k$ neighborhood of $J_0$
with an open set in the Hilbert space $L^2_k(\OTS)$. By the
implicit function theorem, any integrable complex structure
$J=J_0+\mu\in \J^{int}$ with $||\mu||_{L^2_k}$ small is in the
$\G$ orbit of an integrable complex structure $J_0+\nu$ with
$\nu\in Ker Q_0^*$ and $||\nu||_{L^2_k}$ small.
Since both sides of (\ref{lojainfi}) are invariant under the action of $\G$, it suffices to prove it for $\mu\in Ker Q_0^*$.\\

We still need to fix another gauge so that the problem becomes
elliptic. Recall that $\J^{int}$ is the subvariety of $\J$ cut out
by the equation:
$$N(\mu)=\pb_0\mu+[\mu, \mu]=0.$$ We would like to linearize this space to
$Ker\pb_0$.  Let $W=Ker Q_0^*\cap Ker \pb_0$. Consider the
operator
$$\Phi: (W\cap L^2_k(\OTS))\times (Im \pb_0\cap L^2_{k+1}(\OOTS))\rightarrow Im
\pb_0\cap L^2_{k-1}(\OOTS)$$ by sending $(\mu, \alpha)$ to the
orthogonal projection to $Im \pb_0$ of $N(\mu+\pb_0^*\alpha)$.
Since the linearization $$D\Phi_0(\nu, \beta)=\pb_0\pb_0^*\beta
$$ whose second component is an isomorphism, by the implicit
function theorem, for any $\nu\in W\cap L^2_k(\OTS) $ with
$||\nu||_{L^2_k}$ small, there exists a unique
$\alpha=\alpha(\nu)\in Im \pb_0\subset L^2_{k+1}(\OOTS)$ with
$||\alpha||_{L^2_{k+1}}$ small such that $\mu=\nu+\pb_0^*\alpha$
satisfies $\Phi(\mu)=0$. Furthermore, we have
$$||\alpha(\nu)||_{L^2_{k+1}}\leq C\cdot||\nu||_{L^2_k}^2.$$
Define a map $L$ from $B_{\epsilon_1}(W\cap L^2_k(\OTS))$ to
$KerQ_0^*\cap L^2_k(\OTS)$ by sending $\nu$ to $\mu,$ then $L $ is
real analytic and a neighborhood  of $J_0$ in $\J^{int}\cap Ker
Q_0^*\cap L^2_k(\OTS)$ is contained in the image of $L$. Moreover
we have that for all $\nu\in B_{\epsilon_1}W\cap L^2_k(\OTS)$ and
$\lambda\in W\cap L^2_l(\OTS)$(for any $l\leq k$),
\begin{equation}
c_l\cdot||\lambda||_{L^2_l}\leq||(DL)_{\nu}(\lambda)||_{L^2_l}\leq
C_l\cdot||\lambda||_{L^2_l},\label{eq2}
\end{equation}
and
\begin{equation}
c_l\cdot||\lambda||_{L^2_l}\leq||(DL)_{\nu}^*(DL)_{\nu}(\lambda)||_{L^2_l}\leq
C_l\cdot||\lambda||_{L^2_l}. \label{eq3}
\end{equation}
 \\
To be explicit, the differential of $\alpha$ at $\nu$ is given by
$$(D\alpha)_{\nu}(\lambda)=(D\Phi)_{
L(\nu)}(0, -)^{-1}\circ(D\Phi)_{L(\nu)}(\lambda,0).$$ So if we
denote $\mu=L(\nu)$ and $\beta=(D\alpha)_{\nu}(\lambda)$, then
$\beta$ satisfies:
$$\bar{\p}_0\bar{\p}_0^*\beta+\Pi_{Im\pb_0}[\mu, \bar{\p}_0^*\beta]=\pb_0\lambda+\Pi_{Im\pb_0}[\mu, \lambda]=
\Pi_{Im\pb_0}[\mu, \lambda].$$ Thus by ellipticity we obtain for
$\nu$ small that
\begin{equation}\label{eq12}
||(D\alpha)_{\nu}(\lambda)||_{L^2_{l+1}}\leq C\cdot
||\nu||_{L^2_k}\cdot ||\lambda||_{L^2_l}.
\end{equation}
$(\ref{eq2})$ follows from $(\ref{eq12})$ and similarly we can
prove
$(\ref{eq3})$.\\

 Now consider  the
Hilbert space $W\cap L^2_k(\OTS)$ with the constant $L^2$ metric
defined by  $J_0$.  Define the functional $\widetilde{Ca}$ on on a
small neighborhood of the origin in $W\cap L^2_k(\OTS)$ by pulling
back $Ca$ through $L$, i.e.
$$\widetilde{Ca}(\nu)=\frac{1}{2}Ca(L(\nu))=\frac{1}{2}\int (S(L(\nu))-\underline{S})^2\omega^n.$$
It is easy to see that
$$\delta_{\lambda}S(L(\nu))=2Im
\D_{L(\nu)}^*((DL)_{\nu}(\lambda))$$ So the gradient is
$$\nabla \widetilde{Ca}=(DL)_{\nu }^{*}(J\D_{L(\nu)} S(L(\nu))). $$
 We first prove that in a neighborhood of $0$ in $W$,
\begin{equation}||\nabla \widetilde{Ca}(\nu)||_{L^2}\geq C\cdot (\widetilde{Ca}(\nu))^{\alpha}. \label{eq4}
\end{equation}
The linearization of the gradient is the Hessian:
$$H_0:=\delta_{\cdot}\nabla \widetilde{Ca}: L^2_k(W)\rightarrow L^2_{k-4}(W); \lambda\mapsto 2 J_0\D_0\D_0^*\lambda. $$
$H_0$ is an elliptic operator, so it has a finite dimensional
kernel $W_0$ consisting of smooth elements, and $W$ has the
following decomposition:
$$W=W_0\oplus W',$$
where $H_0$ restricts to invertible operators from $L^2_k(W')$ to
$L^2_{k-4}(W')$. So there exists a $c>0$, such that for any
$\mu'\in W'$, we have
$$||H_0(\mu')||_{L^2_{k-4}}\geq
C\cdot||\mu'||_{L^2_k}.$$ By the implicit function theorem, for
any $\mu_0\in W_0$ with $||\mu_0||_{L^2}$\footnote{Since $W_0$ is
finite dimensional, any two norms on it are equivalent. We use the
$L^2$ norm for our later purpose.} small, there exists a unique
element $\mu'=G(\mu_0)\in W'$ with $||\mu'||_{L^2_k}$ small, such
that $\nabla \widetilde{Ca}(\mu_0+\mu')\in W_0$. Moreover the map
$G: B_{\epsilon_1}W_0\rightarrow B_{\epsilon_2}W'$ is real
analytic. Now consider the function
$$f: W_0\rightarrow \R; \mu_0\mapsto \widetilde{Ca}(\mu_0+G(\mu_0)).$$
By construction, this is a real analytic function. For any
$\mu_0\in W_0$, it is easy to
see that $\nabla f(\mu_0)=\nabla \widetilde{Ca}(\mu_0+G(\mu_0))\in W_0$.\\

Now we shall estimate the two sides of inequality $(\ref{eq4})$
separately. For any $\mu\in W$ with $||\mu||_{L^2_k}\leq
\epsilon$, we can write $\mu=\mu_0+G(\mu_0)+\mu'$, where $\mu_0\in
W_0$, $\mu'\in W'$, and
$$||\mu_0||_{L^2_k}\leq c\cdot ||\mu||_{L^2_k},$$
$$||G(\mu_0)||_{L^2_k}\leq c\cdot ||\mu||_{L^2_k},$$
$$||\mu'||_{L^2_k}\leq c\cdot ||\mu||_{L^2_k}.$$
For the left hand side of $(\ref{eq4})$, we have:
\begin{eqnarray*}
\nabla \widetilde{Ca}(\mu)&=&\nabla \widetilde{Ca}(\mu_0+G(\mu_0)+\mu')\\
              &=&\nabla \widetilde{Ca}(\mu_0+G(\mu_0))+\int_0^1
              \delta_{\mu'}\nabla \widetilde{Ca}(\mu_0+G(\mu_0)+s\mu')ds\\
              &=&\nabla f(\mu_0)+\delta_{\mu'}\nabla \widetilde{Ca}(0)+\int_0^1
              (\delta_{\mu'}\nabla \widetilde{Ca}(\mu_0+G(\mu_0)+s\mu')-\delta_{\mu'}\nabla
              \widetilde{Ca}(0))ds\\
\end{eqnarray*}
The first two terms are $L^2$ orthogonal to each other. For the
second term we have
$$||\delta_{\mu'}\nabla \widetilde{Ca}(0)||_{L^2}^2=||H_0(\mu')||_{L^2}^2\geq C\cdot||\mu'||^2_{L^2_4}.$$ For
the last term, we have
$$||\delta_{\mu'}\nabla \widetilde{Ca}(\mu_0+G(\mu_0)+s\mu')-\delta_{\mu'}\nabla
              \widetilde{Ca}(0)||\leq C\cdot||\mu||_{L^2_k}||\mu'||_{L^2_4}\leq C\cdot\epsilon\cdot ||\mu'||_{L^2_4}. $$
Therefore, we have \begin{equation}||\nabla
\widetilde{Ca}(\mu)||_{L^2}^2\geq |\nabla
f(\mu_0)|_{L^2}^2+C\cdot||\mu'||_{L^2_4}^2.\label{eq5}\end{equation}

For the right hand side of $(\ref{eq4})$, we have
\begin{eqnarray*}
\widetilde{Ca}(\mu)&=& \widetilde{Ca}(\mu_0+G(\mu_0)+\mu')\\
       &=& \widetilde{Ca}(\mu_0+G(\mu_0))+\int_0^1 \nabla
       \widetilde{Ca}(\mu_0+G(\mu_0)+s\mu')\mu'ds\\
       &=& f(\mu_0)+\nabla f(\mu_0)\mu'+\int_0^1\int_0^1
       \delta_{\mu'}\nabla \widetilde{Ca}(\mu_0+G(\mu_0)+st\mu')\mu'dtds\\
       &=& f(\mu_0)+H_0(\mu')\mu'+\int_0^1\int_0^1
       (\delta_{\mu'}\nabla \widetilde{Ca}(\mu_0+G(\mu_0)+st\mu')-\delta_{\mu'}\nabla \widetilde{Ca}(0))\mu'dtds
\end{eqnarray*}
So
\begin{equation}
\widetilde{Ca}(\mu)\leq |f(\mu_0)|_{L^2}+
C\cdot||\mu'||_{L^2_4}^2. \label{eq6}
\end{equation}
Now we apply the  Lojasiewicz inequality to $f$, and obtain that
$$|\nabla f(\mu_0)|_{L^2}\geq C\cdot|f(\mu_0)|^{\alpha}, $$
for some $\alpha\in [\frac{1}{2},1)$. Together with $(\ref{eq5})$
and
$(\ref{eq6})$ we have proved $(\ref{eq4})$.\\

To prove $(\ref{lojainfi})$, we need to compare  $||\nabla
Ca(L(\nu))||_{L^2}$ and $||\nabla \widetilde{Ca}(\nu)||_{L^2}$,
i.e. we want
\begin{equation}
||(DL)^{*}_{\nu  }(\D_{L(\nu)}S(L(\nu)))||_{L^2}\leq C\cdot
||\D_{L(\nu)}S(L(\nu))||_{L^2}. \label{eq7}
\end{equation}
We can take $L^2$ decomposition
$$\D_{L(\nu)}S(L(\nu))=(DL)_{\nu}\lambda+\beta,
$$ where $\lambda\in W$ and $\beta\in Ker (DL)^*_{\nu}$. So we just need to prove
$$||(DL)^*_{\nu}(DL)_{\nu}\lambda||_{L^2}\leq C\cdot||(DL)_{\nu}\lambda||_{L^2}$$ for any $\lambda$.
This follows from $(\ref{eq2})$ and $(\ref{eq3})$. $\square$\\

Now we follow the Lojasiewicz arguments. Suppose we have a Calabi
flow $J(t)$ along an integral leaf staying in a $L^2_k$
neighborhood of $J_0$, then by (\ref{Calabi flow J})
$$\frac{d}{dt}Ca(J)^{1-\alpha}=-(1-\alpha)Ca(J)^{-\alpha}||\nabla
Ca(J)||_{L^2(t)}^2\leq -C\cdot||\nabla Ca(J)||_{L^2(t)}.$$ Thus
\begin{equation}
\label{eq8}\int_0^t ||\dot{J}||_{L^2(s)}ds=\int_0^t ||\nabla
Ca(J(s))||_{L^2(s)}ds\leq {C}\cdot
Ca(J(0))^{1-\alpha}.\end{equation} So we get $L^2$ length estimate
for the Calabi flow in terms of the initial Calabi energy. For
$\gamma$ slightly bigger than $\alpha$, we have for
$\beta=2-\frac{\gamma}{\alpha}<1$,
$$\frac{d}{dt}Ca(J)^{1-\gamma}=-(1-\gamma)Ca(J)^{-\gamma}||\nabla
Ca(J)||_{L^2(t)}^2\leq -C\cdot||\nabla Ca(J)||_{L^2(t)}^{\beta}.$$
So for $\beta\in (2-\frac{1}{\alpha}, 1)$ we have
\begin{equation}\label{eq10}\int_0^t ||\dot{J}(s)||_{L^2(s)}^{\beta}ds=\int_0^t ||\nabla
Ca(J(s))||_{L^2(s)}^{\beta}ds\leq C(\beta)\cdot
Ca(J(0))^{1-(2-\beta)\alpha}.\end{equation} Also we have
polynomial decay:
$$\frac{d}{dt}Ca(t)^{1-2\alpha}\geq C>0,$$
so
\begin{equation}
\label{eq9}Ca(J(t))\leq C\cdot (t+1)^{-\frac{1}{2\alpha-1}}.
\end{equation}

Now we define
$$\U_k^{\delta}=\{J\in C^{k,\lambda}(\J^{int})\mid
||\mu_{J}||_{C^{k,\lambda}}\leq \delta \},$$ where again we
identify $J$ close to $J_0$ with $\mu_J\in\Omega_S^{0,1}(T^{1,0})
$.
 Notice that if $\delta\ll1$, then for any tensor $\xi$, the $C^{k,\lambda}_J$ norms defined by $(J,\omega)$
 are equivalent for any $J\in\U_k^{\delta}$. We omit the subscript $J$ if $J=J_0$. Also for $k$ sufficiently large,
  the Sobolev constant is uniformly bounded in $\U_k^{\delta}$.

\begin{thm}
\label{thm2} Suppose $J_0$ is a cscK metric in $\J^{int}$. Then
there exist $\delta_2>\delta_1>0$, such that for any $J(0)\in
\U_k^{\delta_1}$, the Calabi flow $J(t)(t>0)$ starting from $J(0)$
will stay in $\U_k^{\delta_2}$ all the time.
\end{thm}

\emph{Proof}. Choose $\delta>0$ such that the previous a priori
estimates hold in $\U_k^{\delta}$. If suffices to prove that there
exists $\delta_1<\delta_2<\delta$ such that for any Calabi flow
$J(t)$ with $J(0)\in \U_{k}^{\delta_1}$, if $J(t)\in
\U_{k}^{\delta}$ for $t\in [0,T)$, then $J(T)\in
\U_{k}^{\delta_2}$. By lemma \ref{lem2}, for $t\geq 1$ and $l$, we
have
$$||Rm(J(t))||_{C^{l,\lambda}_t}\leq C(l).$$ Now fix $\beta\in (2-\frac{1}{\alpha},1)$, for any $p$,
 there is an $N(p)$(independent of $t\geq 1$), such that the
 following interpolation inequality holds
$$||\dot{J}(t)||_{L^2_p(t)}\leq C(p)\cdot||\dot{J}(t)||_{L^2(t)}^{\beta}\cdot||\D_JS(J)||
_{L^2_{N(p)}(t)}^{1-\beta}\leq C(p)\cdot
||\dot{J}(t)||^{\beta}_{L^2(t)},
$$
 So by (\ref{eq10}) we have
$$ \int_1^T ||\dot{J}(t)||_{L^2_p(t)}dt
\leq  C(p)\cdot Ca(J(1))^{1-(2-\beta)\alpha}\leq  C(p)\cdot
Ca(J(0))^{1-(2-\beta)\alpha}\leq C(p)\cdot\epsilon(\delta_1).$$
Since the Sobolev constant is uniformly bounded in
$\U_k^{\delta}$, we obtain for any $l$,
$$\int_1^T||\dot{J}(t)||_{C^{l,\lambda}_t}dt\leq C(l)\cdot\epsilon(\delta_1).$$
Therefore,
$$||J(T)-J(1)||_{C^{k,\lambda}}\leq\int_1^T||\dot{J}(t)||_{C^{k,\lambda}}dt\leq
\epsilon(\delta_1).$$  By the finite time stability of the Calabi
flow, we have
$$||J(1)-J_0||_{C^{k,\lambda}}=\epsilon(\delta_1).$$
Thus $$||J(T)-J_0||_{C^{k,\lambda}}\leq \epsilon(\delta_1).$$Now
choose $\delta_2=\frac{\delta}{2}$, and $\epsilon(\delta_1)\leq
\delta_2$, then the theorem is
concluded. $\square$\\

From theorem \ref{thm2}, we know the Calabi flow exists globally
in $C^{k,\lambda}$ and thus by sequence converges to $J_{\infty}$
in $C^{k,\beta}$ for $\beta<\alpha$. Now again by the Lojasiewicz
arguments we see the limit must be unique and the
convergence is in a polynomial rate in $C^{k,\lambda}$.\\

Now we  assume that $J_{\infty}=J_0$ is smooth. Then we can prove
smooth convergence. We first use the ellipticity to obtain a
priori estimates in $\U_k^{\delta}$ for $k\gg 1$.  Any $\mu\in
\U_k^{\delta}$ satisfies the following elliptic system:
\begin{equation}
\left\{
 \begin{array}{lll}
          Im \D_0^*\mu=S(\mu)+O(||\mu||_{L^2_2}^2), \\
          Re\D_0^*\mu=Q_0^*(\mu),\\
          \pb \mu+[\mu, \mu]=0.\\
\end{array}\right.
\end{equation}
So we have the following a priori estimate:
\begin{equation}\label{eq11}
||\mu||_{C^{l+2,\alpha}}\leq
C\cdot(||\mu||_{C^{l,\lambda}}+||S(\mu)||_{C^{l,\lambda}}+||Q_0^*(
\mu)||_{C^{l,\lambda}}).
\end{equation}
From the proof of theorem \ref{thm2}, we know that
$||\mu(t)||_{C^{k,\lambda}}$ and $||S(\mu(t))||_{C^{k,\lambda}}$
are uniformly bounded. Since
$$||Q_0^*(\mu(t))||_{C^{k,\lambda}}\leq \int_t^{\infty}||Q_0^*(\dot{\mu}(s))||_{C^{k,\lambda}_s}ds\leq \epsilon(Ca(J(s)))$$
is bounded, we obtain $||\mu(t)||_{C^{k+2,\alpha}}$  bound, so we
can derive smooth convergence by bootstrapping argument.
 This finishes
the proof of theorem \ref{stabilityinfi}.

Theorem \ref{stabilityinfi} has its own interest.  This yields a purely analytical proof
of an extension of a theorem due to Chen \cite{Ch4} and 
Sz\'ekelyhidi \cite{Sz}. This is inspired by an observation of Tosatti \cite{To}. In particular, we do not require the K\"ahler class to be integral. 

\begin{thm}(\cite{Ch4}) For any $J\in\U$, the Mabuchi functional $E$ on the space of K\"ahler
metrics compatible with $J$ is bounded below, and the lower bound
is achieved by the infimum along the Calabi flow initiating from
$J$. 
\end{thm}

\begin{proof} From the proof of theorem \ref{stabilityinfi} we know the Calabi flow
$J(t)\in \J^{int}$ starting from $J$ converges to a limit
$J_{\infty}$ with estimate
$$Ca(J(t))\leq C\cdot (t+1)^{-\frac{1}{2\alpha-1}}.$$
By lemma \ref{Two Calabi flows}, this is equivalent to the Calabi
flow $\phi(t)$ in the space of K\"ahler metrics compatible with
$J$. Then $$E(\phi(t))=E(\phi(0))-\int_0^t Ca(\phi(s))ds\geq
E(\phi(0))-C\cdot \frac{2\alpha-1}{2\alpha-2}\cdot[1-
(t+1)^{\frac{2\alpha-2}{2\alpha-1}}]\geq -C'.$$ For any other
K\"ahler potential $\phi$, we have by lemma \ref{K energy
convexity 1} that $$E(\phi)\geq E(\phi(t))-\sqrt{Ca(\phi(t))}\cdot
d(\phi, \phi(t)).$$ Since
$$d(\phi, \phi(t))\leq d(\phi, \phi(0))+ d(\phi(0), \phi(t))\leq C+\int_0^t\sqrt{Ca(\phi(s))}ds\leq C\cdot[1+ (t+1)^
{\frac{4\alpha-3}{4\alpha-2}}],$$ we have
$$E(\phi)\geq \liminf_{t\rightarrow\infty}E(\phi(t))-C\cdot (t+1)^{-\frac{1}{4\alpha-2}}\cdot [1+(t+1)^
{\frac{4\alpha-3}{4\alpha-2}}]=\lim_{t\rightarrow\infty}
E(\phi(t))$$ is bounded below. 
\end{proof}

\section{Reduced Calabi flow}

In this section we shall discuss a reduced finite dimensional
problem. The usual Kuranishi method provides a local slice as
follows. Assume $J_0$ is cscK. We have as before the following
elliptic complex:
$$C^{\infty}_0(M;\C)\stackrel{\D_0}{\longrightarrow}T_{J_0}\J=\Omega^{0,1}_S(T^{1,0})\stackrel{\pb_0}{\longrightarrow
}\Omega^{0,2}_S(T^{1,0}).$$ Let
$\square_0=\D_0\D_0^{*}+(\pb_0^*\pb_0)^2$, and $H^1=Ker
\square_0$. Let $G$ be the isotropy group of $J_0$, which is the
group of Hamiltonian isometries of $(M, \omega, J_0)$, with Lie
algebra $\g=Ker \D_0\cap C^{\infty}_0(M;\R)$.  By the classical
Matsushima-Lichnerowicz theorem, $Ker \D_0$ is the
complexification $\g^{\C}$ of $\g$, and so the complexification
$G^{\C}$ of $G$ is a subgroup of the group of holomorphic
transformations of $(M, J_0)$, with Lie algebra $\g^{\C}=Ker
\D_0$. Then the linear $G$ action on $H^1$ extends to an action of
$G^{\C}$. For convenience, we include a proof of the following
standard fact.

\begin{lem}\label{Kuranishi}(Kuranishi)
There exists a neighborhood $B$ of $0$ in $H^1$, and a
$G$-equivariant holomorphic embedding
$$\Phi: B\rightarrow \J,$$
such that:\\
(1). $\Phi(0)=J_0$;\\
(2). If $v_1$ and $v_2$ in $B$ are in the same $G^{\C}$ orbit and
$\Phi(v_1)$ is integrable, then $\Phi(v_2)$ is integrable, and
$\Phi(v_1)$ and $\Phi(v_2)$ are in the same $\G^{\C}$ leaf.
Conversely, if $\Phi(v)$ is integrable and $(d\Phi)_v(u)$ is
tangent to the $\G^{\C}$ leaf at $\Phi(v)$,
 then $u$ is tangent to the $G^{\C}$ orbit at $v$.\\
(3). Any integrable $J$ sufficiently close to $J_0$ lies in the
$\G^{\C}$ leaf of some element in the image of $\Phi$.\\
\end{lem}

\begin{proof} We can identify any $J$ close to $J_0$  with an
element $\mu$ in $\Omega_S^{0,1}(T^{1,0})$, and $J$ is integrable
if and only if $$N(\mu)=\pb_0 \mu+[\mu, \mu]=0.$$ We can first
choose a $G$-equivariant holomorphic embedding $\Psi$ from a ball
$B$ in $\OTS$ into $\J$ with $d\Psi_0=Id$, by using the same
``average trick" as in the proof of lemma \ref{average trick}. Let
$$V=\{\mu\in \OTS|\D_0^*\mu=0\},$$
and
$$U=\{\mu\in\Omega_S^{0,1}(T^{1,0})|N(\mu)=0, \D_0^*\mu=0\}.$$
Denote by $G$ the Green operator for $\square_0$ and $H:
\Omega_S^{0,1}(T^{1,0})\rightarrow H^1$ the orthogonal projection.
Then for any $\mu\in U$, we have
$$\mu=G\square_0\mu+H\mu=-G\pb_0^*\pb_0\pb_0^*[\mu,\mu]+H\mu.$$ Define a
$G$-equivariant map $$F: \Omega_S^{0,1}(T^{1,0})\rightarrow
\Omega_S^{0,1}(T^{1,0}); \mu\mapsto \mu+G\pb_0^*\pb_0\pb_0^*[\mu,
\mu],$$ where both spaces are endowed with the Sobolev $L^2_k$
norm. Its derivative at $0$ is the identity map, so by the
implicit function theorem, there is an inverse holomorphic map
$F^{-1}: V_1(\subset \Omega_S^{0,1}(T^{1,0}))\rightarrow
V_2(\subset\Omega_S^{0,1}(T^{1,0}))$. Let $Q$ be restriction of
$F^{-1}$ on $B=V_1\cap H^1$ and $\Phi$ be the composition
$$\Phi: B\rightarrow \J; v\mapsto \Psi\circ Q(v).$$ Since
$H^1$ consists of smooth elements, the image of $\Phi$ also
consists of smooth elements.\\

Now we check $\Phi$ is the desired map. For any $v\in B$, we have
$$\D_0^*Q(v)=-\D_0^*G\pb_0^*\pb_0\pb_0^*[Q(v), Q(v)]=0,$$ and
$$N(Q(v))=-\pb_0G\pb_0^*\pb_0\pb_0^*[Q(v),Q(v)]+[Q(v), Q(v)]=G(\pb_0^*\pb_0)^2[Q(v),Q(v)]
-H[Q(v), Q(v)].$$ So $N(Q(v))=0$ if and only if $H[Q(v), Q(v)]=0$,
as in \cite{Ku}. Therefore a neighborhood of $0$ in $U$ is an
analytic set contained in the image of $Q$. Since  both $\Psi$ and
$F$ are $G$-equivariant and holomorphic, the first part of $(2)$
is true. Following \cite{Sz}, we define a map $P$ from a
neighborhood of $(J_0, 0)$ in $\J\times C^{\infty}_0(M;\C)$ to
$\J$ as follows. Given $\mu\in \OTS$ representing an element in
$\J$ close to $J_0$, and $\phi=\phi_1+\sqrt{-1}\phi_2\in
C^{\infty}_0(M;\C)$ small. There is a family of Hamiltonian
diffeomorphsms $f_t$ with
$$\dot{f}_t=X_{\phi_1}.$$ Denote $J_1=f_1^*J$. Since $\omega_{\phi}=\omega+\sqrt{-1}dJ_1d\phi_2$
is isotopic to $\omega$ through the path
$\omega_t=(1-t)\omega+t\omega_{\phi_2}$. Then there is a canonical
path of diffeomorphisms $g_t$ such that $g_t^*\omega_t=\omega$.
Now $g_1^*J_1$ is the image under $\Psi$ of an element $\mu_1\in
\OTS$. Then define
$$P(\mu, \phi)=G\D_0^*\mu_1.$$ Then $P$ is a smooth function from
$L^2_k(V)\times L^2_k(M;\C)$ to the orthogonal complement
$L^2_k(A^0)$ of $\g^{\C}$ in $L^2_k(M;\C)$. It is easy to
calculate the derivative of $P$ at $(J_0, 0)$ is
 $$(DP)_{0}(\nu, \psi)=G\D_0^*\nu+G\D_0^*\D_0\psi.$$ The derivative with respect to the second variable is surjective
  with a
finite dimensional kernel $0\times \g^{\C}$. Thus by implicit
function theorem, any integrable complex structure close to $J_0$
lies in the $\G^{\C}$ leaf of an element in $U$, and thus is
contained in the $\G^{\C}$ leaf of the
image of $\Phi$. So (3) is proved.\\

It suffices to prove the last statement in $(2)$. Suppose
$\mu=\Phi(v)$, and $\nu=(d\Phi)_v(u)$ is tangent to the $\G^{\C}$
leaf, i.e $\nu=\D_{\mu}\phi$ for some complex valued function
$\phi$. Then $DP_{(\mu,0)}(0, \phi)=0$. On the other hand, the
kernel of $DP_{(\mu,0)}(0, -)$ has the same dimension as $\dim
\g^{\C}$ if $\mu$ is sufficiently close to zero. Thus, $\phi\in
\g^{\C}$ and $u$ is tangent to the $G^{\C}$ orbit of $v$.
\end{proof}

By \cite{D1} the action of $\G$ on $\J$ has a moment map given by
the scalar curvature functional $\mu=S-\underline{S}:
\J\rightarrow C^{\infty}_0(M;\R)$. The downward gradient flow of
$|\mu|^2$ is just the Calabi flow. Now we reduce this flow to a
finite dimensional flow. Note $G$ as a subgroup of $\G$ acts on
$\J$ with induced moment map
$\bar{\mu}=\Pi_{\g}(S-\underline{S})$. It is the $L^2$ projection
of $\mu$ to $\g$ with respect to the natural volume form. We can
consider the gradient flow of $|\bar{\mu}|^2$, whose equation
reads
\begin{equation} \label{reduced Calabi1}
\frac{d}{dt}J=-\frac{1}{2}J\D_J\bar{\mu}(J).
\end{equation}
 If we have a solution to equation (\ref{reduced Calabi1}) such that $J_t$ is integrable for all $t\in[0,T]$, then
 we can translate it to a flow in $\H$ given by
\begin{equation} \label{reduced Calabi2}
\frac{d}{dt}\phi=\Pi_{f_t^*\g}(S(\phi)-\underline{S}),
\end{equation}
where $f_t$ is the family of diffeomorphism satisfying
$$\frac{d}{dt}f_t=-\frac{1}{2}J_tX_{S(J_t)},$$ and the projection is taken with respect to the volume form of $f_t^*\omega$.
 We will study the relation between this flow  and the Calabi flow
later on. Let us call the flow (\ref{reduced Calabi1}) or
(\ref{reduced Calabi2}) the \emph{reduced Calabi flow}. It is the
gradient flow of the norm squared of the moment map of a finite
dimensional
compact group action. \\

Now we can pull back the K\"ahler structure on $\J$ to $B$,
denoted  by $(\tilde{\Omega}, \tilde{J})$. By the previous lemma,
we know $G$ acts on $(B,\tilde{\Omega}, \tilde{J})$
holomorphically and isometrically, with moment map $\tilde{\mu}$
equal to $\Phi^*\bar{\mu}$. We can then study the reduced Calabi
flow on a finite dimensional ambient space $B$. Let $J$ be an
integrable complex structure $J$ close to $J_0$ such that the
Calabi flow $J(t)$ converges to $J_0$. Suppose $J_0$ is not in the
$\G^{\C}$ leaf of $J$. By property (3) in lemma \ref{Kuranishi},
we can smoothly perturb $J(t)$ to $\bar{J}(t)$ in the $\G^{\C}$
orbit such that $\bar{J}(t)=\Phi(v(t))$ for $v(t)\rightarrow 0\in
B$. Since $\dot{\bar{J}}(t)$ is tangent to the $\G^{\C}$ leaf, by
property (2) in lemma \ref{Kuranishi}, we see that $\dot{v}(t)$ is
tangent to the $G^{\C}$ orbit. So $v$ is de-stabilized by $0$ in
$B$ under the $G^{\C}$ action. By our previous study of the finite
dimensional case, the reduced Calabi flow starting from $v$ exists
for all time and converges to $0$ in the order
$O(t^{-\frac{1}{2}})$, and the corresponding flow $\hat{J}(t)$ in
$G^{\C}/G$ is asymptotic to a rational geodesic ray $\chi$ which
also degenerate $v$ to zero. We can view $\chi$ as a geodesic ray
in $\H$ as well, so  the reduced Calabi flow in
$\H$ is asymptotic to a smooth geodesic ray with the same
degeneration limit. This needs a bit more clarification. 
First of all, for any element $g$ in $G^{\C}$, one can choose a path $g(t)$ in $G^{\C}$ with $g(0)$ equal to identity and $g(1)=g$. Then we have
$$\frac{d}{dt}g(t)\cdot g(t)^{-1}=\xi(t)+\sqrt{-1}\eta(t).$$
We can choose a path $h(t)$ in $G$ with $h(0)$ being identity, such that 
$$\frac{d}{dt}(h(t)g(t))\in\sqrt{-1}\g.$$
This is equivalent to 
$$\frac{d}{dt}h(t)\cdot h(t)^{-1}+h(t)\xi(t)h(t)^{-1}=0.$$
Now we define a map $F$ from an open set in $G^{\C}/G$ to $\H$ as follows. This open set is a geodesic convex open set $\U$ in $G^{\C}/G$ such that $[g].v$ still lies in the previously constructed Kuranishi slice.  Let $v(t)=g(t). v$, and $J(t)=\Phi(v(t))$. 
Then $J(t)$ are all integrable and 
$$\frac{d}{dt}J(t)=-(\D_{J(t)}\xi(t)+J(t)\D_{J(t)}\eta(t)), $$
where $\xi(t)$ and $\eta(t)$ are viewed as functions on $M$ through the inclusion $\g\subset C^{\infty}_0(M;\R)$. Choose an isotopy of Hamiltonian diffeomorphisms $f_t$ such that 
$$\frac{d}{dt} f_t=X_{\xi(t)}. $$
Then $\widetilde{J}(t)=f_t^*J(t)$ satisfies
$$\frac{d}{dt} \widetilde{J}(t)=\widetilde{J}(t)\D_{\widetilde{J}(t)}\widetilde{\eta}(t),  $$
where $\widetilde{\eta}(t)=f_t^*\eta(t)$.  In fact, $\widetilde{J}(t)=\Phi(h(t)g(t).v)$. Then by Lemma \ref{Two Calabi flows}  if we choose an isotopy of diffeomorphisms $k_t$ with 
$$\frac{d}{dt}k_t=-\nabla_{\widetilde{J}(t)} \widetilde{\eta}(t),$$
then 
$$k_t^* \widetilde{J}(t)=J, $$
and $k_t^*\omega=\omega_t=\omega+\sqrt{-1}\p\pb\phi(t)$.  We define $F([g])$ to be $\phi(1)$. Of course we need to show this is well-defined,  it suffices to show the definition is independent of the path chosen in  $G^{\C}/G$. Since $G^{\C}/G$ is always simply connected, we only to show it is invariant under based homotopy.  Fo this,  we choose a two parameter family $g_{s, t}$  in $G^{\C}$ such that $g_{s,0}$ is equal to identity, and $g_{s,1}=g$. Correspondingly we have $h(s, t)$ in $G$ with $h(s,0)$ equal to identity.   Let $\widetilde{g}_{s,t}=h_{s,t}\cdot g_{s,t}$, then we have
$$\frac{\p }{\p t}\widetilde g_{s,t}\cdot \widetilde g_{s,t}^{-1}=\sqrt{-1} \eta(s,t)\in \sqrt{-1}\g.$$
Also we have 
$$\frac{\p}{\p s} \widetilde g_{s, t}\cdot \widetilde g_{s, t}^{-1}=\xi(s, t)+\sqrt{-1} \zeta(s,t)\in \g\oplus \sqrt{-1}\g. $$
So we have the relation
$$\sqrt{-1}\frac{\p}{\p s} \eta(s,t)=\frac{\p}{\p t}\xi(s,t)+\sqrt{-1}\frac{\p}{\p t}\zeta(s,t)+[\sqrt{-1}\eta(s,t), \xi(s,t)+\sqrt{-1}\zeta(s,t)].$$
In particular
$$\frac{\p}{\p s}\eta(s,t)=\frac{\p}{\p t}\zeta(s,t)+[\eta(s,t), \xi(s,t)].$$
Also $\xi(s,0)=\zeta(s,0)=\xi(s,1)=\zeta(s,1)=0$. 
Let $J_{s,t}=\Phi(g_{s,t}.v)$, and $f_{s,t}$ be the two parameter family of diffeomorphisms obtained by fixing $s$ and  integrate along the $t$ direction as before. In particular, $f(s,0)$ is equal to identity for all $s$. 
We compute
$$\frac{\p}{\p s}\frac{\p}{\p t} f_{s,t}^*\omega=-\frac{\p}{\p s} f_{s,t}^* dJ_{s,t} d \eta(s,t)=-\frac{\p}{\p s} dJd f_{s,t}^*\eta(s,t).$$
We have
\begin{eqnarray*}
&&\frac{\p}{\p s} f_{s,t}^* \eta(s,t)\\=&&f_{s,t}^*(\frac{\p}{\p s} \eta(s,t)+\L_{J_{s,t}\nabla_{s,t}{\xi(s,t)}-\nabla_{s,t}{\zeta(s,t)}}\eta(s,t))
\\=&&f_{s,t}^*(\frac{\p}{\p s}\eta(s,t)+\{\xi(s,t), \eta(s,t)\}-\langle \nabla_{s,t}\zeta(s,t), \nabla_{s,t}\eta(s,t)\rangle)
\\=&&f_{s,t}^*(\frac{\p}{\p t}\zeta(s,t)-\langle \nabla_{s,t}\zeta(s,t), \nabla_{s,t}\eta(s,t)\rangle))\\
=&&\frac{\p}{\p t} F_{s,t}^* \zeta(s,t).
\end{eqnarray*}
Thus 
$$\frac{\p}{\p s}|_{t=1} f_{s,t}^*\omega=-dJd(\int_0^1 \frac{\p}{\p t} f_{s,t}^*\zeta(s,t)dt)=-dJd(f_{s,1}^*\zeta(s,1))=0.$$
Thus the map $F$ depends only on the point $[g]$, not on the path chosen. So $F$ is a well-defined smooth map. From this it is clear that $F$ is a local isometric embedding, in particular, the image is totally geodesic. Thus we  have proved that the reduced Calabi flow in $\H$ is asymptotic to a smooth geodesic ray with the same degeneration limit. By Section 4.2 this geodesic ray is indeed rational, i.e. extends to a $\C^*$ action. Then 
it  follows from arguments in \cite{Sz} that
$\chi$ is tamed
by a smooth test configuration, so it is tamed by a bounded geometry in the sense of \cite{Ch3}.\\

To prove that the Calabi flow is asymptotic to the reduced Calabi
flow, we need to  generalize lemma \ref{MGS} to the infinite
dimensional case. Then by the same argument as before, together
with lemma \ref{K energy convexity} that the Mabuchi functional
is weakly convex, one can show

\begin{lem}\label{MGS infi} Let $\hat{J}(t)$ be the reduced Calabi
flow as before and $\hat{\phi}(t)$ be the corresponding flow in
$\H$. Then for any Calabi flow path $\phi(t)\in \H$, we have for
all $t$ that
$$d(\phi(t), \hat{\phi}(t))\leq C.$$
\end{lem}

The proof will be given in the appendix.
Combining all these we arrive at the following theorem:

\begin{thm}\label{parallel infi}
Let $(M,\omega_0, J_0)$ be a csc K\"ahler manifold. Let $J$  be a
complex structure in $\J$ close to  $J_0$ and the Calabi flow
starting from $J$ converges to $J_0$ at the infinity. Suppose
$J_0$ is not in the $\G^{\C}$ leaf of $J$. Then there is a smooth
geodesic ray $\phi(t)$ in the space of K\"ahler metrics
$\H_{\omega, J}$ which is tamed by bounded geometry and
degenerates $J$ to $J_0$ in the space $\J$. Furthermore, $\phi(t)$
is asymptotic to the Calabi flow with respect to the
Mabuchi-Semmes-Donaldson metric in the sense of definition
\ref{parallel curve}.
\end{thm}

\section{Relative Bound for parallel Geodesic rays}
It is well-known that in a Riemannian manifold with non-positive
curvature, the distance between two geodesics is a convex
function. In this section we first justify this property for the
infinite dimensional space $\H$.

\begin{lem}\label{convexity}
Let $\phi_1(t)$ and $\phi_2(t)$ be two $C^{1,1}$ geodesics in
$\H$, then $d(\phi_1(t),\phi_2(t))$ is a convex function of $t$.
\end{lem}

\begin{proof}. We first assume both geodesics are $C^{\infty}$. Let
$\gamma_{\epsilon}(t,s)$ be the $\epsilon$-geodesic connecting
$\gamma_1(t)$ and $\gamma_2(t)$(see \cite{Ch1}), then
\begin{eqnarray*}
\frac{d^2}{dt^2}L(\gamma_{\epsilon}(t))&=& \int_0^1
\frac{1}{|\gamma_{\epsilon, s}|}\{
{|\gamma_{\epsilon,ts}^{\perp}|^2}-R(\gamma_{\epsilon,s},\gamma_{\epsilon,t})\}ds
+ \frac{1}{|\gamma_{\epsilon,s}|}\langle \gamma_{\epsilon,s},
\gamma_{\epsilon, tt}\rangle|^1_0\\ &-& \int_0^1
\frac{\langle\gamma_{\epsilon,ss},
\gamma_{\epsilon,tt}\rangle}{|\gamma_{\epsilon,s}|}+\frac{\langle\gamma_{\epsilon,s},\gamma_{\epsilon,ss}\rangle\langle
\gamma_{\epsilon,s},\gamma_{\epsilon,tt}\rangle}{|\gamma_{\epsilon,s}|^3}ds
\end{eqnarray*}
Along the $\epsilon$-geodesics, we have
$$|\gamma_{\epsilon,ss}|=\sqrt{\int_0^1(\phi_{\epsilon,ss}-\nabla_{\phi_{\epsilon,s}}\phi_{\epsilon,s})
^2\omega_{\phi_{\epsilon}}^n} \leq C(t) \sqrt{\epsilon},$$  where
$C(t)$ is uniformly bounded if $t$ varies in a bounded interval.
Also
$$|\gamma_{\epsilon,tt}|\leq C(t), $$ and
$$|\gamma_{\epsilon,s}|\rightarrow L_t,$$ uniformly for $s\in[0,1]$
and $t$ bounded. Therefore, we have
$$\frac{d^2}{dt^2}L(\gamma_{\epsilon}(t))
\geq -C(t)\sqrt{\epsilon}, $$ so for any $a\leq b$,
$$L_{\epsilon}(ta+(1-t)b)\leq
tL_{\epsilon}(a)+(1-t)L_{\epsilon}(b)+C\sqrt{\epsilon}(t-a)(b-t).$$
Let $\epsilon\rightarrow 0$,
$$L(ta+(1-t)b)\leq
tL(a)+(1-t)L(b).$$ So $L(t)$ is still a convex function, and the
argument of the lemma yields the same conclusion. \\

In the general case we need to define the distance between two
$C^{1,1}$ potentials, which is just the infimum of the length of
all $C^{1,1}$ paths connecting the two points. Clearly the
distance between any two points is always non-negative.\\

Now we assume $\phi_1$ and $\phi_2$ are $C^{1,1}$ but $\phi_i(0)$
and $\phi_i(1)$ are smooth, we want to prove for $t\in[0,1]$,
\begin{equation}
L(t)\leq (1-t)L(0)+tL(1).
\end{equation}
To prove this, choose a $\delta$-geodesic $\phi^i_{\delta}$
approximating $\phi_i$ with endpoints fixed. Let $\phi_{\epsilon,
\delta}(t,s)$ be the geodesic connecting $\phi^1_{\delta}(t)$ and
$\phi^2_{\delta}(t)$, and $L_{\epsilon, \delta}(t)$ be its length.
Then similar calculation  shows that
$$\frac{d^2}{dt^2}L_{\epsilon, \delta}(t)\geq -C\sqrt{\delta} -C(\delta, t)\sqrt{\epsilon},$$
So
$$L_{\epsilon, \delta}(t)\leq (1-t)L_{\epsilon, \delta}(0)+tL_{\epsilon,
\delta}(1)+\frac{1}{2} (C\sqrt{\delta} +C(\delta,
t)\sqrt{\epsilon})t(1-t) .$$ Let $\epsilon\rightarrow0$, we have
$$L_{\delta}(t)\leq
(1-t)L_{\delta}(0)+tL_{\delta}(1)+C\sqrt{\delta} .$$ Let
$\delta\rightarrow0$, we get the desired inequality. So
the theorem is true in this case.\\

If $\phi_i(0)$ and $\phi_i(1)$ are not assumed to be smooth, we
can approximate them weakly in $C^{1,1}$ by smooth potentials
$\phi_i^{\epsilon}(0)$, $\phi_i^{\epsilon}(1)$ respectively. Let
$\phi_i^{\epsilon}(t)$ be the geodesic connecting
$\phi_i^{\epsilon}(0)$ and $\phi_i^{\epsilon}(1)$. Then we know
$d(\phi_1^{\epsilon}(t), \phi_2^{\epsilon}(t))$ is a convex
function. By maximum principle for the Monge-Amp\`ere equations,
we know
$$|\phi_i^{\epsilon}(t)-\phi_i(t)|_{C^0}\leq \max(|\phi_i^{\epsilon}(0)-\phi_i(0)|_{C^0},
|\phi_i^{\epsilon}(1)-\phi_i(1)|_{C^0}).$$ Hence
$|\phi_i^{\epsilon}(t)-\phi_i(t)|_{C^0}\rightarrow0$, in
particular, $d(\phi_i^{\epsilon}(t),\phi_i(t))\rightarrow0$.
Therefore, $d(\phi_1^{\epsilon}(t), \phi_2^{\epsilon}(t))$
converges uniformly to $d(\phi_1(t), \phi_2(t))$. So the latter is
also convex.
\end{proof}

\begin{lem}\label{distance} If $\phi_1$ is in $\H$(i.e. $\phi_1$ is smooth and $\omega_1$ is positive) and $\phi_2$ is $C^{1,1}$, then
$d(\phi_1,\phi_2)=0$ if and only if $\phi_1=\phi_2$.
\end{lem}
\begin{proof} We
can choose $C^{\infty}$ potential $\phi_2^{\epsilon}$ converging
to $\phi_2$ weakly in $C^{1,1}$ as $\epsilon\rightarrow0$. Then by
\cite{Ch1},
$$d(\phi_1,
\phi_2^{\epsilon})\geq \max(\int_{\phi_1\geq
\phi_2^{\epsilon}}(\phi_1-\phi_2^{\epsilon})\omega_{\phi_1}^n,
\int_{\phi_2^{\epsilon}\geq\phi_1}(\phi_2^{\epsilon}-\phi_1)\omega_{\phi_2^{\epsilon}}^n)$$
Let $\epsilon\rightarrow 0$, we get
$$d(\phi_1,\phi_2)\geq \max(\int_{\phi_1\geq\phi_2}(\phi_1-\phi_2)\omega_1^n, \int_{\phi_2\geq\phi_1}
(\phi_2-\phi_1)\omega_2^n ).$$ So if $d(\phi_1,\phi_2)=0$, then
$$\int_{\phi_1\geq\phi_2}(\phi_1-\phi_2)\omega_1^n=0, $$
and $$ \int_{\phi_2\geq\phi_1} (\phi_2-\phi_1)\omega_2^n =0.$$ The
first equation implies $\phi_1\leq \phi_2$. The second equation
implies that $$\int_{\phi_2>\phi_1} \omega_2^n=0.$$ Let
$\Omega=\{x\in M|\phi_2(x)>\phi_1(x)\}$. Then by Stokes' formula,
\begin{eqnarray*}
\int_{\Omega}\omega_1^n&=&\int_{\Omega}\omega_1^n-\omega_2^n\\&=&
\int_{\Omega}\sqrt{-1}\p\pb (\phi_1-\phi_2)\cdot\sum_{j=0}^{n-1}
\omega_1^j\wedge\omega_2^{n-1-j}\\&=&\int_{\p\Omega}\sqrt{-1}\pb(\phi_1-\phi_2)\cdot\sum_{j=0}^{n-1}
\omega_1^j\wedge\omega_2^{n-1-j}\\&=&0.
\end{eqnarray*}
So $\Omega$ is empty. Thus $\phi_1=\phi_2$.
\end{proof}

\begin{cor}\label{C0 bound}
Let $\phi_1$ be a geodesic ray tamed by bounded geometry(see
\cite{Ch3}), and $\phi_2$ another geodesic ray parallel to
$\phi_1$ with $\phi_2(0)$ smooth. Then $\phi_1-\phi_2$ has a
uniform relative $C^{1,1}$ bound(with respect to
$\omega_{\phi_1}$).
\end{cor}

\begin{proof} By \cite{Ch3}, there is a $C^{1,1}$ geodesic ray
$\phi_3$ emanating from $\phi_2(0)$ such that
$|\phi_3(t)-\phi_1(t)|_{C^{1,1}_{\phi_1}}\leq C$. Thus
$d(\phi_2(t), \phi_3(t))$ is uniformly bounded. Since
$\phi_2(0)=\phi_3(0)$, by lemma \ref{convexity}, $d(\phi_2(t),
\phi_3(t))=0$. Lemma \ref{distance} then implies
$\phi_2(t)=\phi_3(t)$. So
$|\phi_2(t)-\phi_1(t)|_{C^{1,1}_{\phi_1}}\leq C$.
\end{proof}

\begin{cor} \label{cor3}Let $\gamma_1(t)$ and $\gamma_2(t)$ be two smooth paths in
$\H$ with $d(\gamma_1(t), \gamma_2(t))$ uniformly bounded. Suppose
$\phi(t)$ is a smooth geodesic ray in $\H$ asymptotic to
$\gamma_1$, then it is also asymptotic to $\gamma_2$.
\end{cor}
\begin{proof} Let $\gamma_i(t,s)$ be the geodesic connecting
$\phi(0)$ and $\gamma_i(t)$ parametrized by arc-length. Fix $s$,
by assumption, $d(\gamma_1(t,s), \phi(s))\rightarrow 0$ as
$t\rightarrow\infty$. So in particular, $d(\phi(0),
\gamma_1(t))\rightarrow\infty$. Suppose $d(\gamma_1(t),
\gamma_2(t))\leq C$. Choose $T$ large enough so that $d(\phi(0),
\gamma_1(T))\gg s+C$. Then $d(\gamma_1(T, T-C), \gamma_2(T,
T-C))\leq 4C$. By lemma \ref{convexity}, as $T\rightarrow\infty$,
$$d(\gamma_1(T, s), \gamma_2(T,s))\leq \frac{s}{T}\cdot
4C\rightarrow 0.$$ By definition, $\phi(t)$ is asymptotic to
$\gamma_2$.
\end{proof}

Similarly we can prove
\begin{cor} \label{cor4}Let $\gamma(t)$ be a smooth path in $\H$ which is
asymptotic to two smooth geodesic rays $\phi_1(t)$ and
$\phi_2(t)$. Then $\phi_1$ and $\phi_2$ are parallel, i.e.
$d(\phi_1(t), \phi_2(t))$ is uniformly bounded. If we assume one
of them is tamed by bounded geometry, say $\phi_1$ then by
corollary \ref{C0 bound},
$|\phi_1(t)-\phi_2(t)|_{C^{1,1}_{\phi_1}}\leq C$.
\end{cor}

\section{Proof of the main theorems}
Now we proceed to prove the main theorems.

\begin{lem} \label{lem3}Suppose $g_i$ is a sequence of Riemmanian metrics on
a manifold $M$. If there are two sequences $f_i$ and $h_i$ of
diffeomorphism of $M$  such that $f_i^*g_i\rightarrow g_1$, and
$h_i^*g_i\rightarrow g_2$ in $C^{\infty}$, then $f_i\circ
h_i^{-1}$ converges by subsequence to a diffeomorphism $f$ in $
C^{\infty}$ with $f^*g_2=g_1$.
\end{lem}

The proof is standard using compactness. We omit it here.

\begin{cor}\label{cor2}
The quotient $\J/\G$ is Hausdorff in the $C^{\infty}$ topology.
\end{cor}

\begin{lem}($C^0$ bound implies no K\"ahler collapsing)\label{compactness}
 Suppose there are two sequences $\phi_i$, $\psi_i\in \H$
converging in the Cheeger-Gromov sense, i.e. there are two
sequences of diffeomorphisms $f_i$, $h_i$ such that
$$f_i^*(J, \omega_{\phi_i})\rightarrow (J_1, \omega_1)$$
and $$h_i^*(J, \omega_{\psi_i})\rightarrow (J_2, \omega_2)$$ in
the $C^{\infty}$ topology. If $|\phi_i-\psi_i|_{C^0}\leq C$, then
$|\phi_i-\psi_i|_{C^k_{\omega_{\phi_i}}}$ is bounded for all $k$,
and there is a subsequence ${k_i}$ such that $f_{k_i}^{-1}\circ
h_{k_i}$ converges in $C^{\infty}$ to a diffeomorphism $f$ with
$f^*J_1=J_2$ and
$f^*\omega_1=\omega_2+\sqrt{-1}\p_{J_2}\pb_{J_2}\phi$.
\end{lem}
The proof is quite standard now, given the volume estimates in
\cite{CH1}. We will omit it here.
\begin{figure}
 \begin{center}
  \psfrag{A}[c][c]{$\chi_1(t)$}
  \psfrag{B}[c][c]{$\phi_1(t)$}
  \psfrag{C}[c][c]{$\phi_2(t)$}
  \psfrag{D}[c][c]{$\chi_2(t)$}
  \includegraphics[width=0.8 \columnwidth]{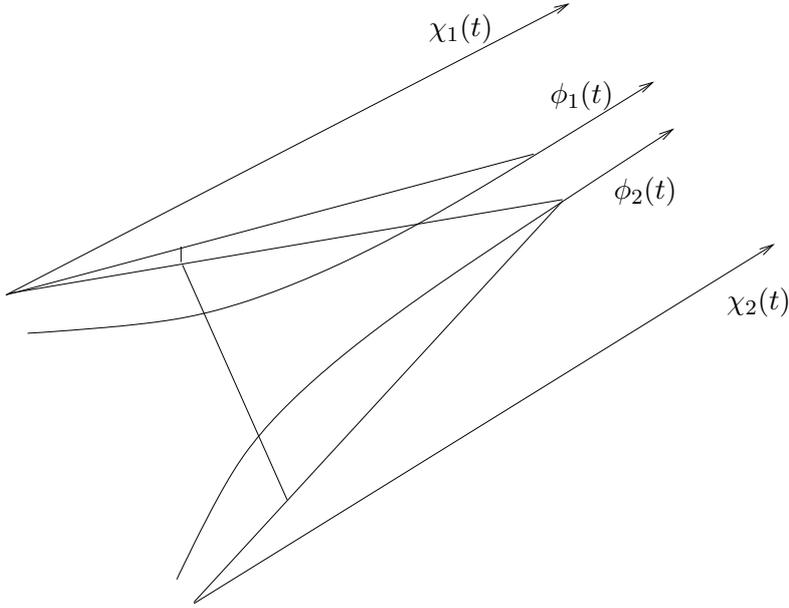}
  \caption{Calabi flows and asymptotic geodesic rays}
  \label{fig: parallel2}
 \end{center}
 \end{figure}

\begin{proof} (of theorem \ref{Theorem 1}). We may assume $J_1$ and $J_2$ are not in the $\G^{\C}$ leaf of $J$,
the proof in the other case is similar. We proceed by
contradiction. Suppose $J_1$ and $J_2$ were not in the same $\G$
orbit. Then by corollary \ref{cor2} we can assume there are
disjoint $\G$ invariant neighborhoods $\U_1$, $\U_2$ of $J_1$,
$J_2$ respectively. Pick $J'_i$ in the intersection of $\U_i$ with
$\G^{\C}$ leaf of $J$. Now by theorem \ref{stabilityinfi}, we know
that the Calabi flow $J_i(t)$ starting from $J'_i$ exits globally
and converges to $J_i(\infty)\in \U_i$. So $J_1(\infty)$ and
$J_2(\infty)$ are not in the same $\G$-orbit either. By theorem
\ref{parallel infi}, the corresponding Calabi flow $\phi_i(t)$ in
the space of K\"ahler metrics is asymptotic to a smooth geodesic
ray  which also degenerates some other $\hat{J}_i$ to
$J_i(\infty)$. Since $J'_1$ and $J'_2$ are both in the $\G^{\C}$
leaf of $J$, we can pull everything back to $J$ and then we have
two Calabi flows $\phi_i(t)$ each asymptotic to a smooth geodesic
ray $\chi_i(t)$ tamed by bounded geometry. By \cite{CC},
$d(\phi_1(t), \phi_2(t))$ is decreasing, so by corollary
\ref{cor3},  $\phi_1(t)$ is also asymptotic to $\chi_2(t)$. By
corollary \ref{C0 bound} and corollary \ref{cor4}
$$|\chi_1(t)-\chi_2(t)|_{C^{1,1}}\leq C.$$ So lemma
\ref{compactness} implies, there is no K\"ahler collapsing, and
there is a diffeomorphism $f$ with $f^*J_1(\infty)=J_2(\infty)$,
and $f^*\omega=\omega+\sqrt{-1}\p\pb \phi$.
Since$(f^*\omega,J_2(\infty))$ and $(\omega, J_2(\infty)$ are both
csc K\"ahler structures in the same K\"ahler class, by theorem
\ref{CT}, there is a diffeomorphism $h$ with
$h^*J_2(\infty)=J_2(\infty)$ and $h^*f^*\omega=\omega$, so
$(f\circ h)^*(\omega,J_1(\infty))=(\omega, J_2(\infty))$. Contradiction.\\
\end{proof}

\begin{proof} (of theorem \ref{Theorem 2}). Suppose $f_i^*(\omega_{\phi_i}, J)\rightarrow (\omega_1, J_1)$, and
$h_i^*(\omega_{\psi_i}, J)\rightarrow (\omega_2, J_2)$. Since
$[\omega]$ is integral, we see that
$[f_i^*\omega_{\phi_i}]=[\omega_1]$ for $i$ large enough, so we
can further assume that $f_i^*\omega_{\phi_i}=\omega_1$, and
$h_i^*\omega_{\psi_i}=\omega_2$. Then we can follow the proof of
theorem \ref{Theorem 1}.
 \end{proof}

\begin{proof} (of corollary \ref{Corollary 1}). Suppose $f_i^*J\rightarrow
J_1$. Since $c_1(J_1)>0$, we have $c_1(f_i^*J)>0$, and we can
choose a sequence of K\"ahler metrics $\omega_i$ in $c_1(J)$ such
that $f_i^*\omega_i\rightarrow \omega_1$. Then we can apply
theorem \ref{Theorem 2}.
\end{proof}

\section{Further Discussions}
There are also some further interesting questions.
\begin{prob} A general notion of optimal degenerations and its
relation to the Calabi flow. Generalize the theorem to the
uniqueness of some ``canonical" objects in the closure, allowing
the occurrence of  singularities. On the other hand, by the Yau-Tian-Donaldson conjecture, one would like to know if there is a direct algebraic-geometric counterpart of Theorem \ref{Theorem 1}, i.e. whether a K-polystable adjacent K\"ahler structure is unique.
\end{prob}

\begin{prob} Quantization approach(\cite{D4}, \cite{Fi}).  In the case of discrete automorphism group, Donaldson \cite{D4} proved the existence of cscK metric implies asymptotic Chow stability. Theorem \ref{CT} in this case follows immediately. It looks like one can use the finite dimensional Kempf-Ness theorem to deal with Theorem \ref{Theorem 1} also. However, this can not be straightforward. The reason is that for an adjacent cscK K\"ahler structure whose underlying complex structure is different from the original one, the automorphism group can not be finite; and it is known that 
the existence of cscK metric(or even KE metric)  does not necessarily imply asymptotic
Chow poly-stability, see the recent counter-example in \cite{OSY}, \cite{DZ}. It seems to the authors that more delicate work is required to proceed by the quantization method.
\end{prob}

\begin{prob} It follows from our result that Tian's conjecture in \cite{Ti2} is
 likely to hold for cscK metrics(The original conjecture allows mild singularities). In the case of
general extremal metrics, we might need to modify the statement in
Tian's conjecture a bit. This can be easily seen in the
corresponding finite dimensional analogue. In that case any
gradient flow can be reversed and we can get critical points in
the limit along both directions of the flow. Clearly they are not
in the same $G$ orbit and therefore ``adjacent" critical point is
not necessarily unique. In our infinite dimensional case, the
naive uniqueness also fails for adjacent extremal metrics.  Such
examples were already implicit in Calabi's seminal paper
\cite{Ca2}. Namely, we consider the blown up of $\P^2$ at three
distinct points $p_1$, $p_2$ and $p_3$(denoted by $Bl_{p_1, p_2,
p_3}(\P^2)$), then by \cite{APS}, the class
$\pi^*[\omega_{FS}]-\epsilon^2([E_1]+[E_2]+[E_3])$ contains
extremal metrics for $\epsilon$ small enough. If $p_1$, $p_2$ and
$p_3$ are in general position(i.e. they do not lie on a line),
then $Bl_{p_1, p_2, p_3}(\P^2)$ are all bi-holomorphic and by
\cite{Ca2} the classes
$\pi^*[\omega_{FS}]-\epsilon^2([E_1]+[E_2]+[E_3])$ have vanishing
Futaki invariant thus the extremal metrics are cscK. If $p_1$,
$p_2$ and $p_3$ lie on a line, Calabi pointed out in \cite{Ca2}
that there is no cscK metric due to  the Lichn\'erowicz-Matsushima
theorem. It is easy to see that for a fixed K\"ahler class
$\pi^*[\omega_{FS}]-\epsilon^2([E_1]+[E_2]+[E_3])$, the extremal
metrics in the case $p_1$, $p_2$ and $p_3$ lie on a line are
adjacent to the cscK metrics in the case $p_1$, $p_2$, $p_3$ are
in general position. So we can find proper extremal metrics even
adjacent to cscK metrics. C. Lebrun also pointed out to us another
example, where we can look at the Hirzebruch surfaces $F_{2n}$ of
even degree. If $n>m$, then with appropriate polarization,
$F_{2n}$  is adjacent to $F_{2m}$, while in \cite{Ca1}, Calabi
explicitly constructed  extremal metrics in any K\"ahler
classes.\\

The problem where the uniqueness fails can be seen from the fact
that our proof depends on the Calabi flow in an essential way.
Since the Calabi flow can only detect de-stabilizing extremal
metrics, we might want to consider only the uniqueness of
de-stabilizing(i.e. energy minimizing) extremal metrics, as a
modification of Tian's conjecture. This idea of de-stabilizing
extremal metrics has already been implicitly discussed in
\cite{Ch3}.
\end{prob}

\begin{prob}
The integrality assumption in theorem \ref{Theorem 2} is just for
fixing the symplectic form. It seems possible to remove this
assumption.
\end{prob}

\appendix
\section{Marle-Guillemin-Sternberg normal form}

 In this
appendix, we shall give a proof of the Marle-Guillemin-Sternberg
normal form theorem for a Hamiltonian group action in the finite
dimensional case(lemma \ref{MGS}). We shall also consider an
infinite dimensional case for our purpose(lemma \ref{MGS infi}).
We suppose that there is a compatible complex structure, which in
general we can not standardize
without some ``errors".\\

\subsection{Model case}

We first look at a prototype. Suppose $\omega$ is a K\"ahler
metric defined in a neighborhood of $0$ in $\C^n$. Then we can not
trivialize both the complex structure and the symplectic structure
simultaneously, however, we can make either of them
standard, with appropriate control on the other.\\\\
 First, it is easy
choose a holomorphic coordinate such that
$$\omega=\omega_0+O(|z|^2),$$ where $\omega_0$ is the standard
symplectic form on $\C^n$. In this way the complex structure is
made standard, while the error on the symplectic form is
quadratic.\\\\
Now we denote $\alpha=\omega-\omega_0.$ Let $f_t:z\rightarrow tz$
be the contraction map. Then
$$\alpha=f_1^*\alpha-f_0^*\alpha=d\theta,$$
where $\theta=\int_0^1f_t^*(X\lrcorner \alpha)dt$, and
$X=z\frac{\p}{\p z}+\bar{z}\frac{\p}{\p \bar{z}}$. So
$$\theta=O(|z|^3).$$
Let $\omega_t=(1-t)\omega+t\omega_0$, then
$\phi_t^*\omega_t=\omega_0$, where $\phi_t$ is the isotopy
generated by the vector fields $Y_t$ satisfying
$$Y_t\lrcorner\omega_t=-\theta.$$ Thus, $Y_t=O(|z|^3)$ and so
$$\phi_t(z)=z+O(|z|^3),$$  and $$\phi_t^*J_0-J_0=O(|z|^2).$$ In this way the symplectic structure is standard, with an
quadratic error on the complex structure.\\

\subsection{Proof of lemma 4.5}
 Suppose a compact group $G$ acts  on a K\"ahler manifold
 $(M, \Omega, J)$ with moment map $\mu$, and
$z_0$ is a zero of the moment map, but not fixed by the whole
group $G$. We denote by $G_0$ the isotropy group of $z_0$ and
$\g_0$ its Lie algebra. We also fix an $\text{Ad}_{G}$-invariant
metric on $\g$. Now consider $Y=G\times_{G_0}(\m\oplus N)$. Here
$N$ is the orthogonal complement of $\g.z_0$ in
$(\g.z_0)^{\omega_0}$, and $\m$ is the orthogonal complement of
$\g_0$ in $\g$. We identify $\rho\in\m$ with $\tilde{\rho}\in
J_0\cdot (\g.z_0)$ through
\begin{equation} \label{Identification}
\langle\rho, \eta\rangle=\Omega_{0}(\tilde{\rho}, X_{\eta}).
\end{equation}
This also induces an identification between $\m$ and $\g.z_0$
which is different from the one coming from the action.
 $G_0$ acts on $(N, \Omega_N=\Omega_0|_N)$ linearly with a
natural moment map $\mu_N$. $Y$ is in fact the symplectic quotient
of $G\times (\g_0\oplus\m\oplus N)\simeq T^*G\times N$ by $G_0$.
The induced symplectic form on $Y$ is given explicitly by(see
\cite{OR})
\begin{eqnarray*}
&&\tilde{\Omega}_{[g,\rho,v]}((L_g\xi_1,\rho_1,v_1),(L_g\xi_2,\rho_2,v_2))\\\\&&:=
\langle\rho_2+d_v\mu_N(v_2),\xi_1\rangle-\langle\rho_1+d_v\mu_N(v_1),
\xi_2\rangle+\langle\rho+\mu_N(v),[\xi_1,\xi_2]\rangle\\\\&&\ \ \
\ +\Omega_0(X_{\xi_1},
X_{\xi_2})+\Omega_0(v_1,v_2)\\\\&&=\langle\rho_2+d_v\mu_N(v_2),\xi_1\rangle-\langle\rho_1+d_v\mu_N(v_1),\xi_2
\rangle+\langle\rho,[\xi_1,\xi_2]\rangle
+\Omega_{0}(v_1,v_2)+\langle\mu_N(v),
[\xi_1,\xi_2]\rangle\\\\&&=\Omega_0(X_1,X_2)+\langle\rho,[\xi_1,\xi_2]\rangle+(
\langle d_v\mu_N(v_2),\xi_1\rangle-\langle d_v\mu_N(v_1),\xi_2
\rangle ) +\langle\mu_N(v), [\xi_1,\xi_2]\rangle,
\end{eqnarray*}
where we identify $T_g G$ with $\g$ through left translation, and
$X_i=X_{\xi_i}+\alpha_i+v_i$ is viewed as a tangent vector at
$z_0$. The $G$ action on $Y$ is Hamiltonian with moment map:
$$\tilde{\mu}: Y\rightarrow \g; [g, \rho, v]\rightarrow
Ad_{g}^*(\mu_N(v)+\rho).$$ To prove lemma \ref{MGS}, we  need to
trace the proof of the relative Darboux theorem. Since $\Omega$ is
K\"ahler, we can choose holomorphic coordinates on a neighborhood
$V$ of $z_0$ such that $\Omega-\Omega_0=O(r^2)$. Let $exp_{z_0}$
be the exponential map with respect to the metric induced from $J$
and $\Omega$. Then we have
$$exp_{z_0}(\rho+v)=z_0+\rho+v+O(r^3).$$
Consider the map
$${exp}: G\times_{G_0} (\m\oplus N) \rightarrow M; (\xi, \rho, v)\mapsto e^{\xi}.exp_{z_0}(\rho+v).$$
 This is a diffeomorphism from a $G$-invariant neighborhood $U$ of $G\times 0$ to a neighborhood $V$ of $G.z_0$.
  Indeed, its derivative at $[e,0,0]$ is given by
 $$d{exp}_{z_0}: \m\oplus\m\oplus N\rightarrow T_{z_0}M=\m\oplus\m\oplus N;(\xi, \rho, v)\mapsto(L(\xi),\rho, v),$$
 where we have made use of the identification (\ref{Identification}), and $L: \m\rightarrow \m$ is the the automorphism such that
$$(L(\xi), \eta)=g_0(X_{\xi}, X_{\eta})$$ for any $
\xi, \eta\in\m$.
 Denote $\Omega'=exp^*\Omega$ and $J'=exp^*J$, then we have
$$J'_{(0,0,0)}(\xi, \rho, v)=(L^{-1}(\rho), -L(\xi), J_0\cdot v).$$
We can extend $J'$  to an almost complex structure
$\tilde{J}$ defined on $Y$. \\

On $V$, denote by $(z,\bar{z})$ the coordinates for $N$, $x$ for
$\g.z_0$ and $y$ for $W=J_0\cdot (\g.z_0)$. The tangent space at
$z_0$ is naturally identified with $V$. Let $(\frac{\p}{\p v},
\frac{\p}{\p \bar{v}})$, and $\frac{\p}{\p \rho}$ be the vector
fields on $U$ corresponding to $\frac{\p}{\p z}$, $\frac{\p}{\p
\bar{z}}$ and $\frac{\p}{\p y}$(on $\m\oplus N$) respectively and
$\frac{\p}{\p \xi}$ the vector fields induced by left translation
of $\frac{\p}{\p x}\in T_{z_0}V$. These vector fields could also
be viewed as vector fields on $V$ through the map $exp$. Then at
$[e, \rho, v]$ we have
$$\frac{\p}{\p v}=\frac{\p}{\p z}+O(r^2);$$
$$\frac{\p}{\p \bar{v}}=\frac{\p}{\p \bar{z}}+O(r^2);$$
$$\frac{\p}{\p \rho}=\frac{\p}{\p y}+O(r^2);$$
$$L\frac{\p}{\p\xi}=\frac{\p}{\p x}+\xi.y+\xi.z+O(r^2).$$
Now it is easy to see that
\begin{eqnarray*}
\tilde{\Omega}(\frac{\p}{\p z}, \frac{\p}{\p
\bar{z}})&=&\Omega'(\frac{\p}{\p z}, \frac{\p}{\p
\bar{z}})+O(r^2),
\end{eqnarray*}

$$\tilde{\Omega}(\frac{\p}{\p z}, \frac{\p}{\p
z})=\Omega'(\frac{\p}{\p z}, \frac{\p}{\p z})+O(r^2)=O(r^2);$$
$$\tilde{\Omega}(\frac{\p}{\p z}, \frac{\p}{\p
y})=\Omega'(\frac{\p}{\p z}, \frac{\p}{\p y})+O(r^2)=O(r^2);$$
$$\tilde{\Omega}(\frac{\p}{\p y}, \frac{\p}{\p y})=\Omega'(\frac{\p}{\p
y}, \frac{\p}{\p y})+O(r^2)=O(r^2).$$

and \begin{eqnarray*} \tilde{\Omega}(\frac{\p}{\p z}, \frac{\p}{\p
x})&=& \tilde{\Omega}(\frac{\p}{\p v}+O(r^2),
L\frac{\p}{\p\xi}-\xi.y-\xi.z+O(r^2))\\&=& \Omega'(\frac{\p}{\p
z},\frac{\p}{\p x})+O(r);
\end{eqnarray*}
similarly, $$\tilde{\Omega}(\frac{\p}{\p y}, \frac{\p}{\p x})=
\Omega'(\frac{\p}{\p y},\frac{\p}{\p x})+O(r);$$

 \begin{eqnarray*} \tilde{\Omega}(\frac{\p}{\p x}, \frac{\p}{\p x})&=&
\tilde{\Omega}(L\frac{\p}{\p\xi}-\xi.y-\xi.z+O(r^2),
L\frac{\p}{\p\xi}-\xi.y-\xi.z+O(r^2))\\&=&\Omega'(\frac{\p}{\p x},
\frac{\p}{\p x})+O(r).\\
\end{eqnarray*}
Therefore, we obtain:
$$\alpha=\Omega'-\tilde{\Omega}=O(r^2)(dzd\bar{z}+dzdy+d\bar{z}dy+dydy)+O(r)(dzdx+d\bar{z}dx+dydx+dxdx).$$
Now let $f_t:(g, \rho, v)\rightarrow (g, t\rho, tv)$, then
$$X_t=\dot{f_t}=t\rho\frac{\p}{\p \rho}+tv\frac{\p}{\p v}+t\bar{v}\frac{\p}{\p \bar{v}}=ty\frac{\p}{\p y}
+tz\frac{\p}{\p z}+t\bar{z}\frac{\p}{\p \bar{z}}+O(r^2).$$ We have
$$\alpha=d\theta, $$
with \begin{eqnarray*} \theta&=&\int_0^1
f_t^*(X_t\lrcorner\alpha)dt\\\\&=& \int_0^1(ty\frac{\p}{\p y}
+tz\frac{\p}{\p z}+t\bar{z}\frac{\p}{\p \bar{z}}+O(r^2))\lrcorner
[O(r^2)(dzd\bar{z}+dzdy+d\bar{z}dy+dydy)\\&&+O(r)(dzdx+d\bar{z}dx+dydx+dxdx)]dt\\\\&=&O(r^2)dx+O(r^3),
\end{eqnarray*}where the estimate is valid at $[e, \rho, v]$.
Let $\Omega_t=(1-t)\tilde{\Omega}+t\Omega'$, then
$$\phi_t^*\Omega_t=\tilde{\Omega}, $$ where
$\dot{\phi_t}=Y_t$ satisfies
$$Y_t\lrcorner\Omega_t=\theta. $$
Since
$$\Omega_t=\Omega_0+O(r^2)(dzd\bar{z}+dzdy+d\bar{z}dy+dydy)+O(r)(dzdx+d\bar{z}dx+dydx+dxdx).$$
So at $[e, \rho, v]$, we have $$Y_t=O(r^2)\frac{\p}{\p
y}+O(r^3)=O(r^2)\frac{\p}{\p \rho}+O(r^3).$$ Since $Y_t$ is
$G$-invariant, this is also true at $[g, \rho, v]$ for $g$ close
to $Id$. Thus the integral curve of $Y_t$ satisfies
$$v_t=v_0+O(r_0^3);$$
$$\rho_t=\rho_0+O(r_0^2).$$

Therefore,
$$(\phi_t^*J')
\frac{\p}{\p v}={\phi_t^{-1}}_*J'((\phi_t)_*\frac{\p}{\p
v})={\phi_t^{-1}}_*J'(\frac{\p}{\p
v}+O(r^2))={\phi_t^{-1}}_*J(\frac{\p}{\p
z}+O(r^2))=\tilde{J}\frac{\p}{\p v}+O(r^2),$$ and similarly
$$(\phi_t^*J')
\frac{\p}{\p \bar{v}}=\tilde{J}\frac{\p}{\p\bar{ z}}+O(r^2).$$ Let
$\Phi=\phi_1$, then $\Phi^*\Omega'=\tilde{\Omega}$. We get the
required estimate that
$$\Phi^*J'-\tilde{J}=O(r),$$ and
$$\Phi^*J'\cdot X-\tilde{J}\cdot X=O(r^2)|X|,$$ for $X\in N$.
Hence lemma \ref{MGS} is proved.

\subsection{Proof of lemma 6.2}
Now we proceed to our infinite dimensional problem, following the
same route as in the finite dimensional setting. However, there
are a few more technical issues, as we shall see below.
 Suppose $(M,\omega,J_0)$ is a
csc K\"ahler manifold. Then the relevant group $\G$ is the group
 of Hamiltonian diffeomorphisms of $(M,\omega)$, which acts on
the space $\J$ of almost complex structures compatible with
$\omega$. Here in order to apply the implicit function theorem, we
shall put $C^{\infty}$ topology on these infinite dimensional
objects which makes them into tame Fr\'echet spaces(\cite{Ha}).
$\J$ inherits a natural (weak) K\"ahler structure $(\Omega, I)$
from the original K\"ahler manifold $M$. The action of $\G$
preserves the K\"ahler structure and has a moment map given by the
Hermitian scalar curvature functional $m(J)=S(J)-\underline{S}$.
Denote by $G$ the identity component of the holomorphic isometry
group of $(M,\omega,J_0)$. Let $\g$ and $\g_0$ the Lie algebra of
$\G$, $G$ respectively. Then we have an $L^2$ orthogonal
decomposition
$$\g=\g_0\oplus\m,$$ where $\m$ is the image of $Q^*=Re \D^*$.
 We want to show that a neighborhood $V$
of $J_0$ in $\J$ is $\G$-equivariantly Hamiltonian diffeomorphic
to a neighborhood $U$ in
$$Y=\G\times_{G}(\m\oplus N), $$ where
 $\G$ acts adjointly on $\g$ by $$f.\phi=f^*\phi.$$ $N$ is the
orthogonal complement of the image of $\D$ in $\OT$, and $G$ acts
on $N$ by pulling back: $g.\mu=g^*\mu$. This action is Hamiltonian
with moment map given by
$$m_N: N\rightarrow \g_0; (m_N(v),\xi)=\frac{1}{2}\Omega(\xi.v,v). $$
Similar to the finite dimensional case we can define a (weak)
symplectic form on $U$. The left $G$ action on $Y$ is Hamiltonian
with moment map given by
$$\tilde{m}:[g,\rho, v]=g^*(\rho+m_N(v)).$$
The exponential map $\Psi$ on $\J$ with respect to the natural
Riemannian metric is well defined by fiber-wise exponential map of
the symmetric space $Sp(2n)/U(n)$, and it is easy to see that it
is a local tame embedding of a neighborhood of the origin in
$\OTS$ into $\J$. Using the local holomorphic coordinate chart of
$\J$, the K\"ahler form satisfies
$$\Omega_{
\mu}=\Omega_0+O(|\mu|^2).$$ It is also clear that
$$\Psi(\mu)=\mu+O(|\mu|^3).$$
Here the norms on both sides could be taken to be the same. Now we
can define a map
$$\Phi: U\rightarrow \J; [g, \rho, v]\mapsto g^*\Psi(\rho+v).$$

\begin{lem}
$\G$ is a smooth tame Lie group.
\end{lem}

\begin{proof} We first prove it is a smooth tame space. We can
identify a Hamiltonian diffeomorphism $H$ with an exact Lagrangian
graph $G_H$ in $\M=M\times M$, i.e.
$$G_H=\{(x, H(x))|x\in M\}.$$ Here $\M$ is endowed with a canonical symplectic form
$\omega'=\pi_1^*\omega-\pi_2^*\omega$, where $\pi_i$ is the
projection map to the $i$-th factor. A Lagrangian graph is called
\emph{exact} if it can be deformed by exact Lagrangian isotopies
to the identity. We can construct local charts for $\G$ as
follows. Given any $H\in \G$, by Weinstein's Lagrangian
neighborhood theorem(\cite{We}), we can choose a symplectic
diffeomorphism between a tubular neighborhood $\U$ of $G_H$ in
$\M$ and a tubular neighborhood $\V$ of $0$ section in the
cotangent bundle $T^*M$. Then locally any Hamiltonian
diffeomorphism close to $H$ is represented by the graph of an
exact one-form, i.e. the differential of some real valued function
on $M$. So locally $\U$ can be identified with an open subset of
$C^{\infty}_0(M;\R)$. Thus $\G$ is modelled on
$C^{\infty}_0(M;\R)$. Now we check the transition function is
smooth tame.
 In our case locally between any two charts there is a symplectic diffeomorphism of
the cotangent bundle $F:T^*M\rightarrow T^*M$ which is identity on
the zero section. Then the induced transition map is smooth tame,
by observing that the $C^k$ distance between the graph of exact
one-forms $d\phi_1$ and $d\phi_2$ is equivalent to the $C^{k+1}$
distance between $\phi_1$ and $\phi_2$. Similary we can prove that
the group multiplication and inverses are both smooth tame.
 \end{proof}

Since the finite dimensional group $G$ acts smooth tame and freely
on $\G\times(\m\oplus N)$, we know that
$$Y=\G\times_{G}(\m\oplus N)$$ is a tame space with a smooth tame
$\G$- action.

\begin{lem} The $\G$-equivariant map $$\Phi:\G\times_G(\m\oplus N)\rightarrow \J; [g, \rho, v]\mapsto g^*\Psi(\rho+v)$$
is smooth tame with a local smooth tame inverse around $[Id,
0,0]$.
\end{lem}
\begin{proof} It is clear by definition that the map is smooth and
tame. The $k$-th derivative of $\Phi$ is tame of degree $k+1$. To
apply Hamilton's implicit function theorem, we need to study the
derivative of $\Phi$ near $[Id,0,0]$. At $\delta=[g, \rho, v]$, we
denote $\mu=\Phi(\delta)$. Then we have
\begin{eqnarray*}
D_{\delta}\Phi:&& \m\oplus\m\oplus N \rightarrow \OTS; [\phi,
\psi, u]\mapsto
\\&&(Id-\bar{\mu})\circ(I-\mu\circ\bar{\mu})^{-1}\circ
[Q_{\mu}\phi+g^*D\Psi|_{\rho+v}(\sqrt{-1}\D_0\psi+u)]\circ
(Id-\mu)^{-1}.
\end{eqnarray*}
To find the inverse to $D_{\delta}$, we need to first decompose
$\OTS$ into the direct sum of $D\Psi^{-1}\circ Im Q_{\mu}|_{\m}$
and $Ker Q_0^*$ with estimate.
 This can be done using elliptic theory.
We can obtain that
$$\nu=(D\Psi)^{-1}\circ Q_{\mu}\phi+\sqrt{-1}Q_0\psi+\eta,$$where $\eta\in Ker
D_0^*$. Take the map $P_{\mu}: \nu\mapsto (\phi,
\sqrt{-1}Q_0\psi+\eta)$. Then it is smooth tame again by elliptic
estimates. Since the inverse of $D_{\delta}\Phi$ is the
combination of $P$ with some other smooth tame operator, it is
also smooth tame. Then we can apply the Nash-Moser implicit
function theorem(\cite{Ha}) to conclude the lemma.
\end{proof}

As in the finite dimensional case, there is a canonically defined
(weak) symplectic form on $\U$ given by
\begin{eqnarray*}
&&\tilde{\Omega}_{[g,\rho, v]}((L_g\xi_1, \rho_1,v_1), (L_g\xi_2,
\rho_2,v_2))
\\\\&&:=\langle\rho_2+d_v\mu_N(v_2),\xi_1\rangle-\langle\rho_1+d_v\mu_N(v_1),
\xi_2\rangle+\langle\rho+\mu_N(v),
[\xi_1,\xi_2]\rangle+\Omega_0(v_1,v_2)\\\\&&=(\D_0^*\D_0\xi_1,
\rho_2)-(\D_0^*\D_0\rho_1+d_v\mu_N(v_1)-[\D_0^*\D_0\rho+\mu_N(v),\xi_1],
\xi_2)+\Omega_0(v_1,v_2)+(\xi_1, d_v\mu_N(v_2))
\end{eqnarray*}
By the above lemma we can pull back the symplectic form $\Omega$
and the complex structure $I$ to $\U$,  denoted by $\Omega'$ and
$I'$ respectively. There is also a canonical almost complex
structure $I_0$ on $\U$ defined by
$$I_0: \m\oplus\m\oplus N\rightarrow \m\oplus\m\oplus N; (\xi, \rho, v)\rightarrow ((\D_0^*\D_0)^{-1}\rho, -\D_0^*\D_0\xi,
I(0)(v)).$$ It is easy to see that $I'=I_0$ at $[Id, 0,0]$.

\begin{prop} There are  neighborhoods $\U_i, \V_i(i=1,2)$($\U_2\subset \U_1$) of $[Id, 0,0]$ in $Y$ and two $
\G$-equivariant smooth tame maps $$\Sigma_1: \U_1\rightarrow
\V_1,$$ $$\Sigma_2: \V_2\rightarrow \U_2,$$ which fixes the
$\G$-orbit of $[Id, 0,0]$ such that  $\Sigma_1\circ \Sigma_2$
equal to the identity and such that
$$\Sigma_1^*\tilde{\Omega}=\Omega',$$
$$\Sigma_2^*\Omega'=\tilde{\Omega},$$ and for any $X\in N$, and
$[g, \rho, v]\in V_2$,
$$(D\Sigma_1)\circ I'\circ(D\Sigma_2)(X)-I_0(X)=O(r^2)\cdot |X|,$$
and at $[Id, 0,0]$, $$(D\Sigma_1) \circ I'\circ (D\Sigma_2)=I_0.$$
Here the estimate is only in the tame sense, i.e. the norm on the
left hand side might be weaker than that on the right, $r$ is the
norm of $[g, \rho, v]$.
\end{prop}

\begin{proof}  The idea of the proof is the same as the finite dimensional case. The main difficulty
is to show the existence of solutions to the involved O.D.E's in
infinite dimension. Once this is established, then everything else
will follow formally.  First we have ($\mu=\Phi([Id, \rho,v])$)
\begin{eqnarray*}&&
\Omega'_{[g, \rho,v]}((L_g\xi_1,\rho_1,v_1), (L_g\xi_2,\rho_2,
v_2))
\\\\&&=\Omega_{\mu}(\D_{\mu}\xi_1+d\Phi_*(\sqrt{-1}\D_0\rho_1+v_1),
\D_{\mu}\xi_2+d\Phi_*(\sqrt{-1}\D_0\rho_2+v_2))
\\\\&&=-Im(\D_{\mu}\xi_1+d\Phi_*(\sqrt{-1}\D_0\rho_1+v_1),
\D_{\mu}\xi_2+d\Phi_*(\sqrt{-1}\D_0\rho_2+v_2))_{L^2}
\\\\&&=(-Im\D_{\mu}^*\D_{\mu}\xi_1-Re \D_{\mu}^*\circ d\Phi_*(\D_0\rho_1-\sqrt{-1}v_1),
\xi_2)\\\\&&\ \ +(Re\D_0^*\circ (d\Phi_*)^t(\D_{\mu}\xi_1)+Im
\D_0^*\circ(d\Phi_*)^td\Phi_*(\D_0\rho_1-\sqrt{-1}v_1),
\rho_2)\\\\&&\ \ +\Omega_0(
(d\Phi_*)^t(\D_{\mu}\xi_1+d\Phi_*(\sqrt{-1}\D_0\rho_1+v_1)), v_2)
\end{eqnarray*}
As in the finite dimensional case, we need to solve an O.D.E. Let
$\Omega_t=(1-t)\tilde{\Omega}+t\Omega'$. The isotopies
$f_t:[g,\rho, v]\rightarrow [g, t\rho, tv]$ gives rise to
time-dependent vector field $X_t(f_t([g, \rho, v]))=[0, \rho, v]$.
We first need to solve another time-dependent vector field $Y_t$
through the following relation:
\begin{equation}
\Omega_{t\ [g, \rho, v]}(Y_t, Z)=\int_0^1
(\Omega'-\tilde{\Omega})_{[g, s\rho, sv]}((0,\rho, v), {f_s}_*Z)ds
\end{equation}
 Notice that $Y_t$ is $\G$-invariant. So we
can assume $g=Id$. Let $Y_t=(\xi_1, \rho_1, v_1)$ and $Z=(\xi_2,
\rho_2, v_2)$. By choosing $Z$
 arbitrarily, we get the following system of equations:
\begin{eqnarray}\label{equation1}
&&-tIm \D_{\mu}^*\D_{\mu}\xi_1-tRe\D_{\mu}^*\circ
d\Phi_*(\D_0\rho_1-\sqrt{-1}v_1)\nonumber\\
\nonumber\\&&-(1-t)\D_0^*\D_0\rho_1-(1-t)
d_v\mu_N(v_1)+(1-t)[\D_0^*\D_0\rho+\mu_N(v),\xi_1]\nonumber\\
\nonumber\\&&
+\int_0^1 Re\D_{\mu_s}^*\circ {d\Phi_{s}}_*(s\D_0\rho+\sqrt{-1}sv)+\D_0^*\D_0(s\rho)-d_v\mu_N(sv) ds\in \g_0\nonumber\\
\end{eqnarray}
\begin{eqnarray}\label{equation2}
&&tRe\D_0^*\circ (d\Phi_*)^t \D_{\mu}\xi_1+tIm \D_0^*\circ (d\Phi_*)^t\circ (d\Phi_*)(\D_0\rho_1-\sqrt{-1}v_1)\nonumber\\
\nonumber\\&&
+(1-t)\D_0^*\D_0\xi_1-\int_0^1 Im \D_0^*\circ(d\Phi_{s})_*^t\circ (d\Phi_s)_*(s\D_0\rho-\sqrt{-1}sv)ds\in \g_0\nonumber\\
\end{eqnarray}
\begin{eqnarray}\label{equation3}
&&t(d\Phi_*)^t\circ (\D_{\mu}\xi_1+d\Phi_*(\sqrt{-1}\D_0\rho_1+v_1))+(1-t)v_1+(1-t)(d_v\mu_N)^*(\xi_1)\nonumber\\ \nonumber\\
&&-\int_0^1({d\Phi_{s}}_*)^t\circ {d\Phi_s}_*(s\sqrt{-1}\D_0\rho+sv)-svds\in Im \D_0\nonumber\\
\end{eqnarray}
Since $\Omega'$ and $\tilde{\Omega}$ are both non-degenerate, this
system admits a (unique) weak solution. Then applying elliptic
regularity, the solutions are smooth. Next we shall prove that
there are two neighborhoods $\NB_1$, $\NB_2$ of $0$
 in $m\times N$, and a smooth tame map $F$ from $\NB_1$ to $C^{\infty}([0,1], \m\times N)$ such that the time $1$ evaluation
 of the image of $F$ is a smooth tame map  from $\NB_1$ to $\NB_2$ and for any $(\rho, v)\in \NB_1$,
\begin{equation}
\left\{
                                                               \begin{array}{ll}
\frac{d}{dt}F_t(\rho,v)=(\rho_1(t),v_1(t)),   & \hbox{ $t\in [0,1]$ ;} \\\\
          F_0(\rho,v)=(\rho,v) .    & \hbox{}\\

                                                               \end{array}
                                                             \right.
\end{equation}\\
To prove this claim, we shall exploit Hamilton's implicit function
theorem again.  Define a map
$$H: C^{\infty}([0,1],\m\times N)\rightarrow (\m\times N)\times C^{
\infty}([0,1], \m\times N)$$ which sends $( \rho(t), v(t))$ to
$(\rho(0), v(0))\times (\dot{\rho}(t)-\rho_1(t),
\dot{v}(t)-v_1(t))$. It is clear that $H$ is a smooth tame map and
$H(0)=0$. We shall show that for $x=(\rho(t), v(t))$ close to
zero, the derivative of $H$ at $x$ is invertible and its inverse
is smooth tame. Let $\delta x=(\tilde{\rho}(t), \tilde{v}(t))$,
then the derivative of $H$ along $\delta x$ is given by
$(\tilde{\rho}(0), \tilde{v}(0))\times
(\dot{\tilde{\rho}}-\delta{\rho_1}(\tilde{\rho}) ,
(\dot{\tilde{v}}-\delta{v_1}(\tilde{v})
 )$. So the invertibilty of $dH$ is equivalent to the solvability of the Cauchy problem of the following linear system
 along $(\rho(t), v(t))$:
\begin{equation}
\left\{
\begin{array}{ll}
\frac{d}{dt}(\alpha,u)=(\delta\rho_1(\alpha),\delta v_1(u))+(\beta, q), & \hbox{$t\in[0,1]$}\\\\
(\alpha(0), u(0))=(\tilde{\rho}(0), \tilde{v}(0)).&\hbox{}\\
\end{array}
\right.
\end{equation}\\

Thus we need to linearize equations (\ref{equation1}) and
(\ref{equation3}). As a result, we get the following
\begin{equation}
\left\{
\begin{array}{lll}
\dot{\alpha}(t)=A_1(\rho(t), v(t))\alpha(t)+A_2(\rho(t), v(t))u(t)+\beta(t), & \hbox{}\\\\
\dot{u}(t)=B_2(\rho(t), v(t))\dot{\alpha}(t)+C_2(\rho(t),
v(t))\alpha(t)+B_0(\rho(t), v(t))u(t)+q(t),
& \hbox{}\\\\
(\alpha(0), u(0))=(\tilde{\rho}(0), \tilde{v}(0)),&\hbox{}\\
\end{array}
\right.
\end{equation}\\
where $A_i$, $B_i$, $C_i$ are pseudo-differential operators of
order $i$ whose coefficients depend on $(\rho(t), v(t))$. Let
$w(t)=u(t)-B_2(\rho(t), v(t))\alpha(t)$. Then the systems of
equations for $(\alpha(t), w(t))$ become symmetric
 hyperbolic, for which the Cauchy problem is always solvable with estimates, see \cite{AG}. From the proof we can check
  that the
 solution depends tamely on $(\rho(t), v(t))$, $(\beta(t), q(t))$ and the initial condition
 $(\tilde{\rho}(0), \tilde{v}(0))$. So  by Hamilton's implicit function theorem $H$ has a local smooth tame inverse.
 Let $F=H^{-1}(-, 0)$ and the claim is then proved.  Now for any $[g, \rho, v]$ close to $[Id, 0,0]$, we obtain a path
 $(\rho(t), v(t))=F_t(\rho, v)$. Then we can solve the O.D.E $\dot{g}(t)=L_{g(t)}\xi_1(t)$,
  where $\xi_1(t)$ is determined by $(\rho(t), v(t))$.  $[g(t), \rho(t), v(t)]$ is then an integral curve of $Y_t$ by
  the $\G$-invariancy. Now we define
 $$\Sigma_2: \V_2\rightarrow \U_2; [g, \rho, v]\mapsto [g_1, F_1(\rho, v)]. $$
Then from the previous arguments we know that $\Sigma_2$ is smooth
tame and fixes $\G.[Id,0,0]$. Moreover,
$\Sigma^*\Omega'=\tilde{\Omega}$. It follows from equations
(\ref{equation1}), (\ref{equation2}), (\ref{equation3}) that we
have a tame estimate
$$|v_1(t)|\leq C\cdot(|\rho(t)|+|v(t)|)^3.$$ Since
$|(v(t),\rho(t))|\leq C\cdot|(v(0),\rho(0))|$, we obtain
$$|v_1(t)-v_1(0)|\leq C \cdot(|\rho(0)|+|v(0)|)^3.$$
By symmetry, we can obtain the map $\Sigma_1$. Then one can check
that the required estimates hold.
\end{proof}

Now to prove lemma \ref{MGS infi}, we just need to apply the
previous proposition to the path $\hat{J}(t)$, and use exactly the
same argument as in the proof of theorem \ref{thmfinite}.

X.X.C. Department of Mathematics, University of
Wisconsin-Madison, 480 Lincoln Drive, Madison, WI 53706, U.S.A./
xxchen@math.wisc.edu;\\
 New address: Department of Mathematics,
Stony Brook University, Stony Brook, NY 11794,
U.S.A./xiu@math.sunysb.edu.\\\\
S.S.  Department of Mathematics, University of
Wisconsin-Madison, 480 Lincoln Drive, Madison, WI 53706, U.S.A./
ssun@math.wisc.edu; \\ New address: Department of Mathematics,
Imperial College, London SW7 2AZ, U.K./s.sun@imperial.ac.uk.

 \vskip3mm


\begin{thebibliography}{0}
\bibitem[AG]{AG} S. Alinhac, P. G\'erard. \emph{Pseudo-differential Operators and the Nash-Moser Theorem}.
 Graduate Studies in
Mathematics, Volume 82. American Mathematical Society, 2007.
\bibitem[APS]{APS} C. Arezzo, F. Pacard, M. Singer. \emph{Extremal metrics on blow
ups},  Duke Math. J. 157 (2011), no. 1, 1--51
\bibitem[BM]{BM} S. Bando, T. Mabuchi. \emph{Uniqueness of Einstein K\"ahler metrics modulo connected group
actions.} Algebraic geometry, Sendai, 1985, Adv. Stud. Pure Math.,
10 (1987), 11-40.
\bibitem[Be]{Be} B.Berndtsson. \emph{Probability measures
related to geodesics in the space of K¨ahler metrics},
arXiv:0907.1806.
\bibitem[Ca1]{Ca1} E. Calabi. \emph{Extremal K\"ahler metric},  Seminar of Differential
Geometry, ed. S. T. Yau, Annals of Mathematics Studies 102,
Princeton University Press (1982), 259-290.
\bibitem[Ca2]{Ca2} E. Calabi. \emph{Extremal K\"ahler metric, II},  Differential Geometry and
Complex Analysis, eds. I. Chavel and H. M. Farkas, Spring Verlag
(1985), 95-114.
 \bibitem[CC]{CC} E. Calabi, X-X. Chen. \emph{Space of K\"{a}hler metrics and Calabi flow}, J. Differential
Geom. 61 (2002), no. 2, 173--193.
\bibitem[Ch1]{Ch1} X-X. Chen. \emph{The space of K\"{a}hler metrics}, J. Differential Geom. 56 (2000), no. 2, 189--234.
\bibitem[Ch2]{Ch2} X-X. Chen. \emph{On the lower bound of the Mabuchi energy and its
application}, Internet. Math. Res. Notices 2000, no. 12, 607-623.
\bibitem[Ch3]{Ch3} X-X. Chen. \emph{Space of K\"ahler metrics III--The greatest lower bound of the Calabi
energy}, Invent. math. 175(2009), 453-680.
\bibitem[Ch4]{Ch4} X-X. Chen. \emph{The space of K\"ahler metrics (IV)---On the lower bound of the K energy}.
Preprint,  arXiv:0809.4081.
\bibitem[CH1]{CH1} X-X. Chen, W-Y. He. \emph{On the Calabi flow}, Amer. J. Math. 130 (2008), no. 2,
539--570.
\bibitem[CH2]{CH2} X-X. Chen, W-Y. He. \emph{The Calabi flow on K\"ahler surface with bounded Sobolev
constant--(I)}, arxiv:math/0710.5159.
\bibitem[CLW]{CLW} X-X. Chen, H. Li, B. Wang. \emph{K\"ahler-Ricci flow with small initial energy.}
 Geom. Funct. Anal. 18 (2009), no. 5, 1525--1563.
\bibitem[CS]{CS} X-X. Chen, S. Sun. \emph{Space of K\"ahler metrics IV-K\"ahler
quantization}, To appear in Metric and Differential Geometry, a
volume in honor of Jeff Cheeger for his 65th birthday.
\bibitem[CT]{CT} X-X. Chen, G. Tian. \emph{Geometry of
K\"{a}hler metrics and foliations by holomorphic discs},
Publications Math\'ematiques de L'IH\'ES, Vol 107, 1-107.
\bibitem[DZ]{DZ} A. Della Vedova, F. Zuddas, \emph{Scalar curvature and asymptotic Chow stability of projective bundles and blowups}, arXiv:1009.5755 .
\bibitem[D1]{D1} S. K. Donaldson. \emph{Remarks on gauge theory, complex geometry and 4-manifold
topology},  Fields Medalists' Lectures, World Sci. Publ.,
Singapore, 1997, 384-403.
\bibitem[D2]{D2} S. K. Donaldson. \emph{Symmetric spaces, K\"{a}hler geometry and Hamiltonian dynamics}, Northern California Symplectic Geometry Seminar, 13--33, Amer. Math. Soc. Transl. Ser. 2, 196, Amer. Math. Soc., Providence, RI,
 1999.
\bibitem[D3]{D3} S. K. Donaldson. \emph{Scalar curvature and projective embeddings},
 I. J. Differential. Geom. 59 (2001), no. 3, 479--522.
\bibitem[D4]{D4} S. K. Donaldson. \emph{Floer homology groups in Yang-Mills theory.} (English summary)
With the assistance of M. Furuta and D. Kotschick. Cambridge
Tracts in Mathematics, 147. Cambridge University Press, Cambridge,
2002.
\bibitem[D5]{D5} S. K. Donaldson. \emph{Conjectures in K\"ahler geometry, Strings and
geometry}, Clay Math. Proc., 3, Amer. Math. Soc., Providence, RI,
(2004), 71-78.
\bibitem[D6]{D6} S. K. Donaldson.
\emph{Constant scalar curvature metrics on toric surfaces}. Geom.
Funct. Anal. 19 (2009), no. 1, 83--136.
\bibitem[DK]{DK} S. K. Donaldson, P. B. Kronheimer. \emph{The geometry of
four-manifolds.} Oxford Mathematical Monographs. Oxford Science
Publications. The Clarendon Press, Oxford University Press, New
York, 1990.
\bibitem[Fi]{Fi} J. Fine. \emph{Calabi flow and projective embeddings
}, arxiv: math/0811.0155.
\bibitem[Fu]{Fu} A. Fujiki. \emph{The moduli spaces and K\"ahler metrics of polarized algebraic varieties.}
 Sugaku Expositions. Sugaku Expositions 5 (1992), no. 2, 173--191.
\bibitem[FS]{FS}A. Fujiki, G. Schumacher, \emph{The moduli space of extremal compact K\"ahler manifolds
and generalized Weil-Petersson metrics.}
 Publ. Res. Inst. Math. Sci. 26 (1990), no. 1, 101--183.
\bibitem[GS]{GS} V. Guillemin, S. Sternberg. \emph{A normal form for the moment
map.} Differential geometric methods in mathematical physics
(Jerusalem, 1982), 161--175, Math. Phys. Stud., 6, Reidel,
Dordrecht, 1984.
\bibitem[Ha]{Ha} R. S. Hamilton. \emph{The inverse
function theorem of Nash and Moser}. Bull. Amer. Math. Soc. (N.S.)
7 (1982), no. 1, 65--222.
\bibitem[Ke]{Ke} G. R. Kempf. \emph{Instability in invariant theory},
Ann. of Math. (2) 108 (1978), no. 2, 299--316.
\bibitem[Ki]{Ki} F. C. Kirwan. \emph{Cohomology of quotients in symplectic and algebraic
geometry}, Mathematical Notes, 31. Princeton University Press,
Princeton, NJ, 1984.
\bibitem[Ku]{Ku} M. Kuranishi. \emph{New proof for the existence of locally complete families of complex structures.}
 1965 Proc. Conf. Complex Analysis (Minneapolis, 1964) pp. 142--154 Springer,
 Berlin.
\bibitem[Le]{Le} E. Lerman. \emph{Gradient flow of the norm squared of a moment
map}, Enseign. Math. (2) 51(2005), no. 1-2, 117-127.
\bibitem[M1]{M1} T. Mabuchi. \emph{Some symplectic geometry on compact K\"{a}hler manifolds I.}
Osaka J. Math. 24 (1987), no. 2, 227--252.
\bibitem[Ne]{Ne} L. Ness. \emph{A Stratification of the Null Cone Via the Moment
Map}, With an appendix by David Mumford. Amer. J. Math. 106
(1984), no. 6, 1281--1329.
\bibitem[OSY]{OSY} H. Ono, Y. Sano, N. Yotsutani, \emph{An example of asymptotically Chow unstable manifolds with constant scalar curvature}, 
arXiv:0906.3836.

\bibitem[OR]{OR} J. Ortega, T. Ratiu. \emph{Momentum maps and Hamiltonian reduction}. Progress in
Mathematics, 222. Birkh\"auser Boston, Inc., Boston, MA, 2004.
\bibitem[Pa]{Pa} S. Paul. \emph{ Hyperdiscriminant polytopes, Chow polytopes, and Mabuchi energy asymptotics
}. arxiv:math/0811.2548.
\bibitem[Ra]{Ra} J. R\aa de. \emph{On the Yang-Mills heat equation in two and three dimensions}.
J. Reine Angew. Math. 431 (1992), 123--163.
\bibitem[RT]{RT} J. Ross, R. Thomas. \emph{A study of the Hilbert-Mumford criterion for the stability of projective varieties.}
 J. Algebraic Geom. 16 (2007), no. 2, 201--255.
\bibitem[Se]{Se} S. Semmes. \emph{Complex Monge-Amp\`{e}re equations and sympletic
 manifolds},
Amer. J. Math, no. 114, 495--550, 1992.
\bibitem[Si]{Si} L. Simon. \emph{Asymptotics for a class of nonlinear evolution equations, with
applications to geometric problems}. Ann. of Math. (2) 118 (1983),
no. 3, 525--571.
\bibitem[SW]{SW} S. Sun, Y-Q. Wang. \emph{On the K\"ahler-Ricci flow near a K\"ahler-Einstein metric},
arxiv:math/1004.2018.
\bibitem[Sz]{Sz} G. Sz\'ekelyhidi. \emph{The K\"ahler-Ricci flow and K-stability}, Amer. J. Math. 132 (2010), 1077--1090.
\bibitem[Th]{Th} R. Thomas. \emph{Notes on GIT and symplectic reduction for bundles and varieties
}, arxiv: math/0512411.
\bibitem[Ti1]{Ti1} G. Tian. \emph{K\"ahler-Einstein metrics with positive scalar curvature}, Invent. Math. 130 (1997),
1-39.
\bibitem[Ti2]{Ti2} G. Tian. \emph{Extremal metrics and
geometric stability. Special issue for S. S. Chern}. Houston J.
Math. 28 (2002), no. 2, 411--432.
\bibitem[To]{To} V. Tosatti. \emph{The K-energy on small deformations of constant saclar curvature
K\"ahler manifolds}. arxiv: math/1010.1859.
\bibitem[TZ1]{TZ1} G. Tian, X-H. Zhu. \emph{A new holomorphic invariant and uniqueness of K\"ahler-Ricci
solitons}. Comment. Math. Helv., 77 (2002), no. 2, 297-325.
\bibitem[TZ2]{TZ2} G. Tian, X-H.Zhu. \emph{Perelman's $\W$-functional and stability of K\"ahler-Ricci
flow}, arxiv: math/0801.3504.
\bibitem[We]{We} A. Weinstein. \emph{Symplectic manifolds and their Lagrangian submanifolds.}
 Advances in Math. 6 (1971), 329-346.
\end{thebibliography}
\end{document}